\RequirePackage{etex}
\documentclass[11pt, a4paper]{article}
\usepackage[utf8]{inputenc}
\usepackage{inputenc}
\usepackage{amsmath}
\usepackage{amssymb}
\usepackage{upgreek}
\usepackage{graphicx}
\usepackage{authblk}
\usepackage{layout}
\usepackage{dsfont}
\usepackage{bbold}
\usepackage[english]{babel}
\usepackage[top=2.2cm, bottom=3.2cm, left=2.2cm, right=2.2cm]{geometry}
\usepackage{multicol}
\usepackage{amsthm}
\usepackage{mathrsfs}
\usepackage{tikz}
\usepackage{tikz-network}
\usepackage[toc,page]{appendix}
\usepackage{hyperref}
\usepackage{subfig}
\usepackage{complexity}
\usepackage[margin=1cm]{caption}
\usepackage{enumitem}
\usepackage{fontawesome}

\newtheorem{theorem}{Theorem}
\newtheorem{corollary}[theorem]{Corollary}

\newtheorem{open}[theorem]{Open question}

\newtheorem{lemma}[theorem]{Lemma}

\newtheorem{remark}[theorem]{Remark}
\usepackage{todonotes}

\usepackage{makecell}
\usepackage{changepage}
\usepackage{diagbox}

\usepackage{parskip}

\newcommand*\samethanks[1][\value{footnote}]{\footnotemark[#1]}

\newcommand{\popogame}{poset positional game}

\newcommand{\popogames}{poset positional games}
\newcommand{\Popogames}{Poset positional games}
\newcommand{\popoproblem}{{\sc MB Poset Positional Game}}
\newcommand{\poset}{P}
\newcommand{\WS}{{\mathcal{F} }}
\newcommand{\Ga}{{\mathcal{G} }}
\newcommand{\strat}{{\mathcal{S} }}
\renewcommand{\emptyset}{\varnothing}

\newcommand{\QBF}{{\sc 3-QBF}}

\newcommand{\OM}{{\mathcal{M} }}
\newcommand{\OB}{{\mathcal{B} }}
\newcommand{\ON}{{\mathcal{N} }}
\newcommand{\OP}{{\mathcal{P} }}
\newcommand{\OA}{{all }}

\tikzstyle{v}=[circle, inner sep=1pt, minimum size =14 pt, line width = 1pt, draw=black, fill=white, text= black]
\tikzstyle{B}=[circle, inner sep=0, minimum size =6 pt, line width = 1pt, draw=black, fill=blue, text= white]
\tikzstyle{R}=[circle, inner sep=0, minimum size =6 pt, line width = 1pt, draw=black, fill=red, text= white]

\title{Poset Positional Games\thanks{This research was partly supported by the ANR project P-GASE (ANR-21-CE48-0001-01)}}
\author[1]{Guillaume Bagan} 
\author[1]{Eric Duch\^{e}ne}
\author[2]{Florian Galliot}
\author[3]{Valentin Gledel}
\author[4]{Mirjana Mikala\v{c}ki\thanks{Partly supported by Provincial Secretariat for Higher Education and Scientific Research, Province of Vojvodina (Grant No.~142-451-2686/2021). Partly supported by Ministry of Science,
Technological Development and Innovation of Republic of Serbia
(Grants 451-03-66/2024-03/200125 \& 451-03-65/2024-03/200125).}} 
\author[1]{Nacim Oijid}
\author[1]{Aline Parreau}
\author[4]{Milo\v{s} Stojakovi\'c\samethanks[2]} 

\affil[1]{Univ Lyon, CNRS, INSA Lyon, UCBL, Centrale Lyon, Univ Lyon 2, LIRIS, UMR5205, F-69622 Villeurbanne, France.}

\affil[2]{Univ. Bordeaux, CNRS, Bordeaux INP, LaBRI, UMR 5800, F-33400 Talence, France}

\affil[3]{Université Savoie Mont Blanc, CNRS UMR5127, LAMA, Chambéry, F-73000, France}

\affil[4]{Department of Mathematics and Informatics, Faculty of Sciences, University of Novi Sad, Serbia}

\date{}

\begin{document}

\maketitle

\begin{abstract}
We propose a generalization of positional games, supplementing them with a restriction on the order in which the elements of the board are allowed to be claimed. We introduce \popogames, which are positional games with an additional structure -- a poset on the elements of the board. Throughout the game play, based on this poset and the set of the board elements that are claimed up to that point, we reduce the set of available moves for the player whose turn it is -- an element of the board can only be claimed if all the smaller elements in the poset are already claimed.

We proceed to analyse these games in more detail, with a prime focus on the most studied convention, the Maker-Breaker games. First we build a general framework around \popogames. Then, we perform a comprehensive study of the complexity of determining the game outcome, conditioned on the structure of the family of winning sets on the one side and the structure of the poset on the other.
\end{abstract}

\section{Introduction}
\subsection{General motivation}
\textbf{Positional games.} Positional games are a class of combinatorial games that have been extensively studied in recent literature -- see books~\cite{beck} and~\cite{positionalgames} for an overview of the field. They include popular recreational games like Tic-Tac-Toe, Hex and Sim. 
Structurally, a positional game is a pair $(X, \WS)$, where $X$ is a (finite) set that we call \emph{the board}, and $\WS \subseteq2^X$ is a collection of sets that we call \emph{the winning sets}. The pair $(X, \WS)$ is referred to as the \emph{hypergraph of the game}. The game is played in the following way: two players alternately claim unclaimed elements of the board, until all the elements are claimed.
Depending on the way the winner is determined, there are several standard conventions: Maker-Maker games, Maker-Breaker games,  Avoider-Enforcer games, etc.

In \emph{Maker-Maker games}, also known as strong making games, players compete to fill up one of the winning sets i.e.~claim all its elements, and whoever does it first wins. If there is no winner by the time all the board elements are claimed, the game is declared a draw. 
Tic-Tac-Toe (``3-in-a-line'') is a notable representative of this type of games and every child knows the game ends in a draw, provided that both players play optimally. It can be generalized to the $n \times n$ board, going under the name ``$n$-in-a-line'', and it is also known to be a drawn game~\cite{beck} for all $n \geq 3$. More generally, Tic-Tac-Toe can be played on the hypercube $[n]^d$, where the winning sets are all the geometric lines of cardinality $n$. The Hales-Jewett Theorem~\cite{Hales1963} states that, for every $n$, there exists a positive number $d$, also called the Hales-Jewett number $HJ(n)$, such that for all $k\geq d$, the game $[n]^k$ is a first player win.

Generally speaking, the strong making games are natural to introduce and study, and many recreational games are of this type. Hence, it is no wonder that numerous questions about them have been asked in the literature. Yet, it may come as a surprise that the majority of those questions remain unanswered. The thing is, in a typical strong making game, each player's goal can be seen as two-fold: they are simultaneously trying to claim a winning set and to block their opponent's attempts at claiming a winning set. This makes most such games notoriously hard to analyse, and hardly any general tools are known for this purpose. As a consequence, there are very few published results about this convention, when compared to Maker-Breaker games.

In \emph{Maker-Breaker games}, we call the players Maker and Breaker, and they have different goals -- Maker wins if she fills up a winning set, while Breaker wins otherwise, i.e.~if he claims at least one element in each winning set.
Note that no draw is possible in this convention. Maker-Breaker games are the most researched convention of positional games, ever since Erd\H{o}s and Selfridge~\cite{erdos} first introduced them in 1973. Looking at Tic-Tac-Toe on the standard $3 \times 3$ board in the Maker-Breaker convention, it is straightforward to convince oneself that Makers wins when playing first.

The main problem on positional games consists in determining the {\em outcome} i.e.~the identity of the player who has a winning strategy (possibly depending on who starts), if there is one, assuming that both players are playing optimally. The study of the complexity of computing the outcome of a given positional game can be traced back to Schaefer~\cite{schaefer}, who was first to prove that Maker-Breaker games are {\sf PSPACE}-complete, even when the winning sets are of size at most $11$. This was later improved by Rahman and Watson~\cite{rahman2021}, requiring only winning sets of size $6$. On the other hand, Galliot et al.~\cite{Florian2} proved that the outcome of any Maker-Breaker game with winning sets of size at most $3$ can be determined in polynomial time. Maker-Maker games are also known to be {\sf PSPACE}-complete, as shown by Byskov in~\cite{byskov2004maker}. 

\textbf{\Popogames.} Let us now take a closer look at a popular recreational game, Connect-4, which has a lot in common with the Tic-Tac-Toe family of positional games. In Connect-4, two players play on a board that is $7$-wide and $6$-high. They move alternately by placing a token of theirs in a column of their choice. Each placed token drops down with ``gravity'', landing on top of the last previously placed token in that column, or heading all the way to the bottom of the column if it is empty. The first player with tokens on four consecutive positions in a line (vertical, horizontal or diagonal) wins, and, if neither of the players manages to do that by the time all the columns are full, then a draw is declared. Although the game is similar to Tic-Tac-Toe at first glance, there is one crucial difference -- the players cannot choose freely from all the empty positions of the board, as at any point a column offers at most one available position (the lowest token-free position).

The outcome of the game is known as a first player win, as shown by Allis in~\cite{allis1988knowledge}, and independently by Allen~\cite{allen1990expert} who later wrote a whole book~\cite{allen2010complete} about this game. Connect-4 was solved in~\cite{tromp2008solving} for some nonstandard board sizes, and for some other in~\cite{tromp2015}.


Recently, Avadhanam and Jena~\cite{restrictedPosGames} studied \emph{Connect-Tac-Toe}, where they introduce a restriction on which unclaimed elements can be claimed, combining the nature of Connect-4 with the positional game of Tic-Tac-Toe.
They also look at the generalized $[n]^d$ Connect-Tac-Toe, played on a hypercube $[n]^d$, relating it to the restricted version of the $[n]^d$ Tic-Tac-Toe. Additionally, they give a lower bound on an analogue of the Hales-Jewett number, $HJ(n)$, in this setting.


In the present paper, we propose a new framework which, in full generality, enables us to combine the move restrictions with positional games (like it is done in Connect-4). Namely, we introduce \emph{\popogames}, which are positional games with an additional structure -- a poset on the elements of the board. Throughout the game play, based on this poset and the set of elements that are claimed up to that point, we reduce the set of available moves for the player whose turn it is -- an element of the board can only be claimed if all the smaller elements (in the poset) are already claimed.

We proceed to analyse these games in more detail, with a prime focus on the most studied convention, the Maker-Breaker games.  After setting out a formal introduction of \popogames, we go on to build a general framework around them. Then we perform a comprehensive study of the complexity of determining the game outcome, conditioned on the structure of the family of winning sets on the one side and the structure of the poset on the other.







\subsection{Framework of \popogames}

\begin{figure}
    \centering
    \begin{tikzpicture}
                \node[v, anchor = center] (u1) at (0,0.75){a};
                \node[v, anchor = center] (u2) at (1,0.75){b};
                \node[v, anchor = center] (u3) at (0.5,1.5){c};
                \node[v, anchor = center] (u4) at (0.5,2.25){d};
                \draw[->] (u1) -- (u3);
                \draw[->] (u2) -- (u3);
                \draw[->] (u3) -- (u4);
                \node[v, anchor = center] (v1) at (4,0){e};
                \node[v, anchor = center] (v2) at (4,0.75){f};
                \node[v, anchor = center] (v3) at (4,1.5){g};
                \node[v, anchor = center] (v4) at (4.5,2.25){h};
                \node[v, anchor = center] (v5) at (3.5,2.25){i};
                \draw[->] (v1) -- (v2);
                \draw[->] (v2) -- (v3);
                \draw[->] (v3) -- (v4);
                \draw[->] (v3) -- (v5);
                \draw[->] (v1) -- (u4);

    \end{tikzpicture}
    \caption{A poset on nine vertices.}
    \label{fig:simpleposet}
\end{figure}
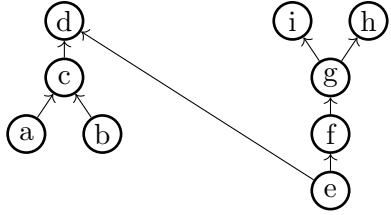

A {\em poset} $P$ on a set $X$ is defined by a partial order relation $\leq$. On top of that, we use $<$ to denote the same relation acting on distinct elements. In this paper, posets will be depicted by directed graphs in the usual way, where the elements are vertices, and two elements $x$ and $y$ satisfy $x\leq y$ if and only if there exists a directed path from $x$ to $y$. For example, in the poset in Figure~\ref{fig:simpleposet}, we have $e\leq i$ and $e\leq d$. Two elements are deemed {\em incomparable} if there exists no directed path between them (e.g.~$b$ and $i$ are incomparable in the poset from Figure~\ref{fig:simpleposet}).

In addition, we will use the following standard definitions about posets: A {\em chain} is a set of elements that are pairwise comparable. In other words, it is a set $S\subseteq X$ such that for any $x,y \in S$, we have $x\leq y$ or $y \leq x$. An {\em antichain} is a set of elements that are pairwise incomparable. In other words, it is a set $S\subseteq X$ such that for any $x,y \in S$, we have $x \not < y$ and $y \not < x$. The {\em height} of a poset is the cardinality of its longest chain. The poset from Figure~\ref{fig:simpleposet} has height $4$. The {\em width} of a poset is the cardinality of its largest antichain. The poset from Figure~\ref{fig:simpleposet} has width $4$ because of the antichain $\{a,b,i,h\}$. Given an element $x$, we say that $y$ is a {\em predecessor} of $x$ (resp. a {\em successor}) if $y<x$ (resp. $x<y$) and there is no other element between $x$ and $y$. An element $x$ is said to be {\em maximal} (resp.~{\em minimal}) in the poset if it has no successor (resp.~no predecessor). 

We are now able to formally define a \popogame~as follows. A {\em \popogame}~is a triple $(X,\WS, \poset)$ where $X$ is a finite set of elements (also called {\em vertices}), $\WS$ is a collection of subsets of $X$ corresponding to the {\em winning sets}, and $\poset$ is a poset on $X$.

The game is played by two players that alternately claim an unclaimed vertex $x$ of $X$ such that all vertices smaller than $x$ have already been claimed. The rest stays the same as before. In the Maker-Maker convention, the first player that fills up a winning set $S\in \WS$ wins. If no player manages to fill up a winning set, the game ends in a draw. In the Maker-Breaker convention, Maker wins if she fills up a winning set at some point during the game, otherwise Breaker wins. During a game, when a player claims a single vertex, we call that a {\em move}, whereas a {\em round} corresponds to a pair of moves made successively by both players. 

Note that a standard positional game 
is a \popogame~where all the vertices are pairwise incomparable. Furthermore, in \popogames, from any given position, the set of moves that are available to the next player forms an antichain of the poset. In particular, there are at most $w$ available moves at any point of a game, where $w$ is the width of the poset.

\begin{remark}\label{rem:simple}
    In the Maker-Breaker convention, we will indifferently use the terms \emph{game} and \emph{game position}. Indeed, any game position that includes moves already played by Maker and Breaker can be simplified to a position without pre-played moves. To do this, it suffices to remove all the vertices already played. The winning sets containing at least one move of Breaker can be removed, and those containing moves of Maker can be reduced to only unclaimed vertices. By definition of the rules, both operations do not affect the winner of the game. With this reduction, all the elements that are available for any player correspond to minimal elements of the poset.
\end{remark}

It is a well-known result that, for any standard positional game played in the Maker-Maker convention, the second player cannot have a winning strategy, so that the only possible outcomes are a first player win $\mathcal{FP}$ or a draw $\mathcal{D}$. Similarly, for any standard positional game played in the Maker-Breaker convention, there are only three possible outcomes (if we do not specify who starts):
\begin{itemize}[noitemsep, topsep = 0cm]
    \item $\mathcal{M}$ if Maker has a winning strategy no matter who starts the game,
    \item $\mathcal{B}$ if Breaker has a winning strategy no matter who starts the game,
    \item $\mathcal{N}$ if the next player (i.e.~the player whose turn it is) has a winning strategy.
\end{itemize}
A simple strategy stealing argument ensures that the second player cannot win in any Maker-Breaker game.

When switching to the framework of \popogames, the property asserting that the second player never wins is not true any more. As a consequence, there are game positions for which the next player may have interest in skipping their turn. The corresponding outcome will be denoted by $\mathcal{P}$ (standing for ``Previous player wins", or equivalently meaning that the second player has a winning strategy). Such phenomena are generally called {\em zugzwangs} in the literature of combinatorial games~\cite{Elkies1999OnNA}. Zugzwangs are known to be hard to deal with, in particular when the games are decomposed into sums (i.e.~disjoint unions) of smaller games, leading to a kind of ``waiting game" within the game. Since \popogames~may have this property, this leads to a more complex analysis of their study from a general point of view. Therefore, it motivates the study of the disjoint union of \popogames~according to the outcome of each game of the union. Figure~\ref{table: unionfirst} gives the possible outcomes in a disjoint union of \popogames, according to the parity of the number of vertices in each game. In particular, one can remark that the outcome $\mathcal{M}$ is often absorbing in a disjoint union of \popogames. The results from this figure will be proved in the appendix (Theorem~\ref{thm:union}).

\begin{figure}[!htb]
\begin{center}
\begin{tabular}{ c || c | c | c | c }
 $even \backslash even$ & $\OM$ & $\ON$ & $\OP$ & $\OB$ \\ 
 \hline \hline
 $\OM$ & $\OM$ & $\OM$        & $\OM$   & $\OM$ \\ \hline  
 $\ON$ & $\OM$ & $\OM, \ON$   & $\OM$   & $\OM, \ON$ \\ \hline 
 $\OP$ & $\OM$ & $\OM$        & $\OP$   & $\OP$ \\ \hline
 $\OB$ & $\OM$ & $\OM, \ON$   & $\OP$   & $\OP, \OB$  \\
\end{tabular}
\end{center}

\vspace{.3cm}
\begin{center}
\begin{tabular}{ c || c | c | c | c }
 $even \backslash odd$ & $\OM$ & $\ON$ & $\OP$ & $\OB$ \\ 
 \hline \hline
 $\OM$ & $\OM$ & $\OM, \ON$         & $\OM$         & $\OM, \ON$ \\  \hline
 $\ON$ & $\OM$ & $\OM, \ON$         & $\OM$         & $\OA$ \\ \hline
 $\OP$ & $\OM$ & $\ON$              & $\OM$         & $\ON$ \\ \hline
 $\OB$ & $\OM$ & $\ON$              & $\OM, \OP$    & $\ON, \OB$  \\
\end{tabular}
\end{center}

\vspace{.3cm}
\begin{center}
\begin{tabular}{ c || c | c | c | c }
 $odd \backslash odd$ & $\OM$ & $\ON$ & $\OP$ & $\OB$ \\ 
 \hline \hline
 $\OM$ & $\OM$      & $\OM$              & $\OM$         & $\OM, \ON$ \\  \hline
 $\ON$ & $\OM$      & $\OM, \OP$         & $\OM, \ON$    & $\OA$ \\ \hline
 $\OP$ & $\OM$      & $\OM, \ON$         & $\OM$         & $\OM, \ON$ \\ \hline
 $\OB$ & $\OM, \ON$ & $\OA$              & $\OM, \ON$    & $\OA$  \\
\end{tabular}
\end{center}
    \caption{ Possible outcomes for a disjoint union of two \popogames~depending on the parity of the number of vertices in both games. We write "$\OA$" when all four outcomes are possible.}
    \label{table: unionfirst}
\end{figure}

  If we have a disjoint union of \popogames, because of the existence of zugzwangs, components without a winning set cannot be directly removed from the union, as they can be considered as waiting moves. Yet, only the parity of such components plays a role in the outcome. 

\begin{lemma}\label{lem:emptyUnion}
    Let $\Ga=\Ga_1\cup \Ga_2$ be a \popogame~that is a disjoint union of two \popogames\ $\Ga_1$ and $\Ga_2$, such that $\Ga_2$ has no winning set. Then the component $\Ga_2$ can be removed from $\Ga$ if its board is even, or replaced by a single vertex if its board is odd, without changing the outcome of the game.
\end{lemma}
\begin{proof}
	Define $\Ga'$ as $\Ga_1$ (if the board of $\Ga_2$ is even) or $\Ga_1$ plus a new isolated vertex $x$ (if the board of $\Ga_2$ is odd), and fix the identity of the first player. We will show that the winning player, i.e. the player who has a winning strategy $\strat$ in $\Ga'$, can also win in $\Ga$ using a strategy that we now describe.
	\begin{itemize}
		\item Firstly, suppose that the board of $\Ga_2$ is even. The winning player can simply follow the strategy $\strat$, except when the opponent plays inside $\Ga_2$, in which case the winning player claims an arbitrary vertex in $\Ga_2$ (as made possible by the parity). Since there is no winning set in $\Ga_2$, this ensures a win in $\Ga$.
		\item Now, suppose that the board of $\Ga_2$ is odd. The strategy is the same, except that we need to handle two events E1 and E2 which may occur at some point:
		\begin{itemize}[noitemsep,nolistsep]
			\item E1: The winning player is instructed by $\strat$ to claim $x$.
			\item E2: The opponent claims the last remaining vertex in $\Ga_2$.
		\end{itemize}
		If E1 occurs before E2, then the winning player can simply claim an arbitrary vertex in $\Ga_2$, after which the parity ensures that E2 never occurs. If E2 occurs before E1, then $x$ is unclaimed in $\Ga'$, so the winning player can imagine that the opponent claimed $x$, after which E1 obviously never occurs. \qedhere
	\end{itemize}
\end{proof}

From now on, unless explicitly stated, we will only consider the Maker-Breaker convention of poset positional games.
Though the computation of the outcome is generally the main issue when investigating such games, this study can, under certain circumstances, be reduced to the case of a particular player starting the game. Indeed, 
starting with a game $\Ga=(X,\WS,\poset)$ where Maker (resp.~Breaker) is the first player and claims a vertex $x$, we get a resulting game $\Ga'=(X',\WS',\poset')$ where Breaker (resp.~Maker) is the first player, defined by $X'=X-x$, $\poset'=\poset-x$ and $\WS'=\{W \setminus \{x\}|\; W \in \WS\}$ (resp.~$\WS'=\{W|\; W \in \WS, x \not\in W\}$). Therefore, when studying a class of games which is stable under Maker's (resp.~Breaker's) moves in terms of that update, we may freely assume that Breaker (resp.~Maker) is the first player, up to considering all possibilities for their opponent's first move otherwise. In this paper, all studied classes will be stable under Breaker's moves, so we will always assume that Maker is the first player. Therefore, the decision problem that will be mainly investigated is the following:
\medskip

\noindent \popoproblem~\\
{\bf Input:} A \popogame{} $\Ga=(X,\WS,\poset)$. \\
{\bf Output:} The player having a winning strategy (i.e.~Maker or Breaker) when Maker starts.
\medskip

Since standard positional games are included in \popogames, the problem \popoproblem\ is {\sf PSPACE}-complete from a result of Schaefer~\cite{schaefer}.

\subsection{Exposition of the results} 

The main objective of this paper is to consider the complexity of \popoproblem~related to some parameters of the instance. More precisely, we have chosen to focus on the properties of the poset, as it is the main distinctiveness of the current contribution, comparing the results we obtain to the previous results about standard positional games. In Section~\ref{sec:height}, we will firstly examine the problem depending on the height of the poset. As it is already known that, even for height $1$ (i.e. all the elements are pairwise incomparable) and winning sets of size $6$, the problem is {\sf PSPACE}-complete \cite{rahman2021}, we will also refine our classification according to the number and/or the size of the winning sets. The main contribution of this section addresses the case where there is only one winning set of size~$1$. By adding such a condition, the problem admits a complexity jump between instances of height $2$, proved to be polynomial, and height $3$, proved to be {\sf NP}-hard.

Section~\ref{sec:width} deals with the width of the poset. As the case of width $1$ is straightforward ($P$ is made of only one chain, so the order of the moves is completely predetermined), we start by considering instances with width $2$. We show that \popoproblem\ is {\sf PSPACE}-hard in this case, even if all the winning sets are of size $3$. This illustrates a major difference with standard positional games, which are known to be tractable in this setting. We then give a polynomial-time algorithm that solves the general case where the width of the poset and the number of winning sets are fixed. 

Section~\ref{sec:disjoint} is devoted to the case where the poset is a union of disjoint chains. This case is a direct generalization of the game Connect-4. Our first result is a full characterization of the outcome when all the winning sets are of size 1. When the winning sets are of size at most $2$, things are more tricky, but we do provide a polynomial-time algorithm when the width of the poset (i.e.~the number of chains in the union) is fixed. 
Finally, by adding the restriction that the height of each chain is at most $2$, we get a polynomial-time algorithm regardless of the number of chains.

Table~\ref{tab:results} summarizes all the results known about the complexity of \popoproblem, referring either to the literature of standard positional games, or to results proved in the current paper. In the table, recall that the complexity class {\sf XP} defines the class of problems parameterized by a parameter $k$ and that can be solved in time $O(|X|^{f(k)})$, where $|X|$ is the size of the instance and $f$ is a computable function. 
In addition to these results, the choice of focusing on the Maker-Breaker convention instead of the Maker-Maker one is reinforced in Appendix by Theorem~\ref{thm:MMsize1depth2}, which yields an instance (height $2$ and all non-minimal elements are winning sets of size $1$) for which the resolution is PSPACE-complete in the Maker-Maker convention, whereas it is in polynomial time in the Maker-Breaker convention (see Theorem~\ref{thm:topverticesWS}).

\begin{table}[htb]
\begin{center}
    \begin{tabular}{!{\vrule width 1.2 pt}c!{\vrule width 1.2 pt}c|c|c|c!{\vrule width 1.2 pt}}
        \noalign{\hrule height 1.2 pt}
        \diagbox{ \\[3mm] winning  \\ sets\\[1mm]}{{} \\[-1mm] poset\;\\[-5mm]} 
        & general & height~$h$& width~$w$ & \makecell{disjoint chains} \\ \noalign{\hrule height 1.2pt}
         general & {\sf PSPACE}-c~\cite{rahman2021} & \makecell{$h=1$:\\ {\sf PSPACE}-c~\cite{rahman2021}} &\makecell{$w=1$: {\sf P} \\[1mm] $w = 2$:\\ {\sf PSPACE}-c \\ (Th.~\ref{thm:psp_w2}) }& \makecell{$h=1$: \\{\sf PSPACE}-c~\cite{rahman2021} } \\ \hline
         size $s$ & \makecell{$s=3$:\\ {\sf PSPACE}-c  \\ (Th.~\ref{thm:psp_w2})}& \makecell{$h=1$, $s=6$:\\{\sf PSPACE}-c~\cite{rahman2021} \\[1mm] $h=2$, $s=3$: \\ $\NP$-h~(Th.~\ref{thm:nph_d2})} & \makecell{$w=2$, $s=3$:\\{\sf PSPACE}-c \\ (Th.~\ref{thm:psp_w2})} & \makecell{$s = 1$: \\{\sf P} (Th.~\ref{thm:disjoint-size1})\\[1mm]$s=2$:\\{\sf XP} by $w$ (Th.~\ref{thm:disjoint-size2})\\[1mm] $h=2$, $s=2$:\\{\sf P} (Th.~\ref{thm:poly_disj_d2})}\\ \hline
         \Gape[0.3cm]{number} $m$ & \makecell{$\NP$-hard\\ (Th.~\ref{thm:nph_d3})} & \makecell{$m=1$, $h=3$:\\$\NP$-hard (Th.~\ref{thm:nph_d3})} & \makecell{ ${\sf XP}$ (Th.~\ref{thm:xp_wl})}& \makecell{$m=1$:\\ {\sf P} (Th.~\ref{thm:disjoint-size1} \\and Remark~\ref{rem:ws1k})}\\ \hline
        \makecell{ number $m$ \\ and size $s$} & \makecell{$m=1$, $s=1$:\\ $\NP$-hard\\ (Th.~\ref{thm:nph_d3})} & \makecell{$h=2$, $m=1$, $s=1$: \\ {\sf P} (Th.~\ref{thm:poly_d2}) \\[1mm] $h=3$, $m=1$, $s=1$: \\ $\NP$-hard (Th.~\ref{thm:nph_d3})} & \makecell{${\sf XP}$ (Th.~\ref{thm:xp_wl})} & \makecell{${\sf XP}$ (Cor.~\ref{cor:disjointboundednumbersize})}\\ \noalign{\hrule height 1.2pt}
    \end{tabular}
    \end{center}
    \caption{Complexity of \popoproblem{} depending on some parameters of the poset and of the collection of winning sets.}
    \label{tab:results}
\end{table}

\section{Posets of small height}\label{sec:height}


\subsection{Posets of height 2} 

\subsubsection{Winning sets of size 1} \label{ss:ws-size-1}

We first look at what happens when all the winning sets are of size 1. This case is rarely straightforward, and often deserves to be studied depending on some parameters of the poset. In addition, it is closely correlated to the case where there is a unique winning set (of any size). The following remark explains the link between both situations.

\begin{remark}\label{rem:ws1k}
    Up to switching the roles of the players, a \popogame~with a single winning set of size $k$ is equivalent to a \popogame~with $k$ winning sets of size $1$.
\end{remark}

In the above remark, the equivalence means that the two games have the same outcome. Indeed, Maker wins when there is a single winning set of size $k$ if and only if she manages to claim all the elements of this set. In the second game, a win of Breaker consists in claiming all the winning sets (of size $1$), thus corresponding to the equivalence. Note that Breaker is starting in one of the two games, but this is not a pitfall by Remark~\ref{rem:simple} since both these classes are stable under Breaker's moves. \\

We continue with a minor result which is useful when dealing with winning sets of size 1. Given a poset $P$ and a vertex $x$, we denote by $p(x)$ the number of vertices that are not greater or equal to $x$
i.e.~$p(x):=|X|-|\{y|\; y \geq x\}|$. This quantity is used in the following lemma, which yields a necessary and sufficient condition for Maker to win when there is a unique winning set of size $1$.

\begin{lemma}\label{lem:predoneWS}
  Let $\Ga=(X,\{\{x\}\},\poset)$ be a \popogame. Suppose that, after some rounds of play in $\Ga$, all predecessors of $x$ have been claimed apart from one, which we denote by $y$. Then, in this position, Maker wins if and only if $p(y)$ is odd.
\end{lemma}

\begin{proof}
  At this point in the game, neither player wants to claim $y$, since it would free up $x$ and the opponent would win on the spot. Therefore, the game is just a parity game, with both players claiming all the vertices they can without claiming $y$. There are exactly $p(y)$ such vertices that can be claimed, including all the vertices that have already been claimed. If $p(y)$ is odd (resp.~even), then Maker (resp.~Breaker) will claim the last vertex that can be claimed without claiming $y$, after which Breaker (resp.~Maker) will be forced to claim $y$ and Maker (resp.~Breaker) will claim $x$ to win.
  \end{proof}

\begin{theorem}\label{thm:poly_d2}
 \popoproblem{} can be solved in polynomial time when: the poset is of height 2, there is only one winning set, and this winning set is of size 1. 
\end{theorem}

\begin{proof}
Let $\Ga=(X,\{\{x\}\},\poset)$ be a \popogame{} with $x\in X$ and $\poset$ has height 2.
If $x$ is a minimal element, then Maker trivially wins. Thus we can assume that $x$ has at least one predecessor. One can partition the predecessors of $x$ into two sets $M$ and $B$ where $M$ (resp.~$B$) contains the predecessors $y$ of $x$ such that $p(y)$ is even (resp.~odd). 

We claim that Maker wins if and only if $|M|\leq |B|$. Indeed, assume that $|M|\leq |B|$. Then Maker repeatedly claims an unclaimed vertex of $M$ until all vertices of $M$ are claimed. At this point $B$ contains at least one unclaimed vertex. Thus, the last predecessor of $x$ to be claimed will be in $B$.
By Lemma~\ref{lem:predoneWS}, this means that Maker wins. If $|M|>|B|$, then similarly Breaker wins by claiming vertices of $B$ -- they will be all claimed before the vertices of $M$ are, and the last predecessor of $x$ to be claimed will be in $M$. We again conclude with Lemma~\ref{lem:predoneWS} that Breaker wins.
\end{proof}

The situation is already much more complicated with two winning sets of size 1, and we leave it as an open problem.

\begin{open}
What is the complexity of \popoproblem{} on a poset of height~2 with two winning sets of size 1?
\end{open}

However, for winning sets of size 1 and a poset of height 2 whose nonminimal elements are all winning sets, one can compute the outcome in polynomial time, as we now show. Note that the same situation is {\sf PSPACE}-complete in the Maker-Maker convention (see Theorem \ref{thm:MMsize1depth2} in Appendix).

\begin{theorem}\label{thm:topverticesWS}
\popoproblem{} can be solved in polynomial time when: the poset is of height 2, all the winning sets are of size 1, and all nonminimal elements are winning sets.
\end{theorem}

\begin{proof}
Let $\poset$ be a poset of height $2$ on a set $X$. We say the minimal elements of $\poset$ are the ``bottom" vertices, and the others are the ``top" vertices. For every relation $x<y$, we know $x$ is a bottom vertex and $y$ is a top vertex. Note that every vertex that is not comparable with any other vertex is a bottom vertex. We say that $x$ is a \textit{private predecessor} of $y$ if $y$ is the only successor of $x$.

The \popogame~we consider is $\Ga=(X,\{\{y\},y $ top vertex$\},\poset)$. We will prove that Maker wins if and only if at least one of the following conditions holds:

\begin{enumerate}
    \item $|X|$ is odd.
    \item There exists a top vertex $y$ which does not have more private predecessors than non-private ones.
\end{enumerate}

First note that these two conditions can be tested in polynomial time.

We prove that Maker wins under any of these two conditions. If condition 1 holds then, since Maker starts, she also plays last. Let $y$ be any top vertex. Maker plays as follows. She never claims $y$ or its last predecessor. When the last predecessor of $y$ is claimed, either it is claimed by Breaker and Maker wins by claiming $y$, or Maker has to claim it. In this case, as the total number of vertices is odd, there are at least two remaining vertices, all of which are top vertices. Therefore, Maker will win with her next move.

  Assume now that Condition 2 holds, and let $y$ be a top vertex which does not have more private predecessors than non-private ones. Maker starts by claiming private predecessors of $y$ until all of them are claimed. Right after the last private predecessor of $y$ is claimed, since there are at least as many non-private predecessors of $y$ by assumption, there will be still some unclaimed predecessors of $y$ that are not private. 
    
    Then, as long as Breaker does not claim the last predecessor of $y$, Maker claims an arbitrary bottom vertex which is not the last predecessor of $y$, if possible. When the last predecessor $x$ of $y$ is claimed, either it is claimed by Breaker and Maker wins by claiming $y$, or Maker has to claim it. The fact that Maker has to claim $x$ means that all the other bottom vertices have been claimed, so all the successors of $x$ are now available. Moreover, since $x$ is not a private predecessor of $y$, it has at least two successors. Therefore, even if Breaker claims one of them, Maker will claim another and win.

\medskip

We now prove that Breaker wins when none of the conditions hold, i.e.~when $|X|$ is even, and the property $(*)$ ``all the top vertices have strictly more private predecessors that non-private ones'' holds.  Note that the status of being private/non-private cannot change during the game.

We give a strategy for Breaker that maintains the property $(*)$ throughout. Whenever Maker claims a bottom vertex $x$, Breaker answers as follows:

    \begin{itemize}
        \item If $x$ was the last predecessor of some top vertex $y$, then Breaker claims $y$. Since $(*)$ was true before Maker's turn and $y$ has only one predecessor $x$, we know $x$ must be a private predecessor of $y$. In particular, there is no other top vertex that has been made available when Maker has claimed $x$. Moreover, $(*)$ is still true after Breaker's move since the predecessors of the other top vertices have not changed. 
        \item Otherwise, if $x$ was a private predecessor of some top vertex $y$ that has at least one non-private predecessor $x'$, then Breaker claims $x'$. Note that this way, $(*)$ still holds for $y$ as well as for all the other top vertices.
        \item Otherwise, if there is a bottom vertex that is not a private predecessor of any top vertex, then Breaker claims an arbitrary such vertex. The property $(*)$ still holds after this move.
        \item Otherwise, every bottom vertex is a private predecessor of some top vertex. Since the total number of vertices is even and it is Breaker's turn, there must be a top vertex $y$ with at least two (private) predecessors. Then Breaker claims an arbitrary predecessor of $y$. The property $(*)$ still holds after this move.
    \end{itemize}
    
    By following this strategy, Breaker thus ensures that the property $(*)$ is satisfied after each of his moves, so Maker will never be able to claim a top vertex. Therefore, Breaker wins.
\end{proof}

From Remark~\ref{rem:ws1k}, note that the above result also resolves \popoproblem\ in polynomial time for posets of height 2 and with a unique winning set built with all the nonminimal elements.

\subsubsection{Winning sets of size 3}

\begin{theorem}\label{thm:nph_d2}
    \popoproblem~is {\sf NP}-hard, even when restricted to instances where the poset has height $2$ and all the winning sets are of size 3.
\end{theorem}

\begin{proof}
    The proof is a reduction from $3$-SAT. Let $\phi$ be a $3$-SAT formula. We build a \popogame~$\Ga = (X, \WS, \poset)$ with $\poset$ of height 2 and winning sets of size 3 as follows:
\begin{itemize}
    \item For any variable $x_i$ of  $\phi$, we add four vertices $u_i, v_i, \overline{v_i}$ and $\widetilde{v_i}$ in $X$.
    \item For any clause $C_j = l_{j_1} \vee l_{j_2} \vee l_{j_3}$ in $\phi$, we add the winning set $S_j$ with, for $1\le k \le 3$, $v_{j_k} \in S_j$ if $l_{j_k} = x_{j_k}$, and $\overline{v_{j_k}} \in S_j$ if $l_{j_k} = \neg x_{j_k}$.
    \item For any variable $x_i$ of $\phi$, we add the relations $u_i < v_i$, $u_i < \overline{v_i}$ and $u_i < \widetilde{v_i}$ in $\poset$.
\end{itemize}

We prove that $\phi$ is satisfiable if and only if Breaker wins the \popogame~played on $\Ga$. Suppose first that $\phi$ is satisfiable and let $\nu$ be a valuation that satisfies $\phi$. We consider the following strategy for Breaker:
\begin{itemize}
    \item Whenever Maker claims a vertex $u_i$, if $x_i$ is set to True in $\nu$ then Breaker claims $v_i$, and if $x_i$ is set to False in $\nu$ then Breaker claims $\overline{v_i}$.
    \item Whenever Maker claims a vertex in $\{v_i, \overline{v_i}, \widetilde{v_i}\}$, Breaker claims the remaining vertex in $\{v_i, \overline{v_i}, \widetilde{v_i}\}$.
\end{itemize}

Note that, since Maker starts, these moves are always possible. For each winning set $S_j$, by assumption, $\nu$ satisfies the corresponding clause $C_j$ by assumption. Thus, there exists some $l_{j_k} \in C_j$ that is evaluated to True. If this literal is $x_{j_k}$, then Breaker has claimed $v_{j_k} \in S_j$, and if it is $\neg x_{j_k}$ then Breaker has claimed $\overline{v_{j_k}} \in S_j$. Thus, Breaker wins.
\smallskip

Conversely, suppose that $\phi$ is not satisfiable, and consider the following strategy for Maker:
\begin{itemize}
    \item As her first move, Maker claims some arbitrary $u_j$.
    \item Whenever Breaker claims a vertex of some pair $\{v_i, \overline{v_i}\}$, Maker responds by claiming the other, if it is still available.
    Otherwise, Maker claims an arbitrary vertex.
\end{itemize}

When all the vertices have been claimed, Maker has claimed at least one vertex in each pair $\{v_i, \overline{v_i}\}$. Consider the valuation $\nu$ such that, for all $i$, $x_i$ is True if and only if Breaker has claimed $v_i$. As Breaker could not have claimed both $v_i$ and $\overline{v_i}$ for any $i$, we know that Maker has claimed all the vertices corresponding to literals that are false under $\nu$. On the other hand, since $\phi$ is not satisfiable by assumption, there exists some clause $C_j$ whose literals are all false under $\nu$. This means all three vertices of $S_j$ have been claimed by Maker, so Maker has won the game.
\end{proof}

\subsection{Posets of height 3} 

We now consider posets of height 3. Using Lemma \ref{lem:predoneWS}, we prove that \popoproblem{} is \NP-hard even if there is only one winning set and that winning set has size 1.

\begin{theorem}\label{thm:nph_d3}
\popoproblem\ is {\sf NP}-hard even when restricted to instances where: the poset has height 3, there is only one winning set, and that winning set is of size 1.
\end{theorem}

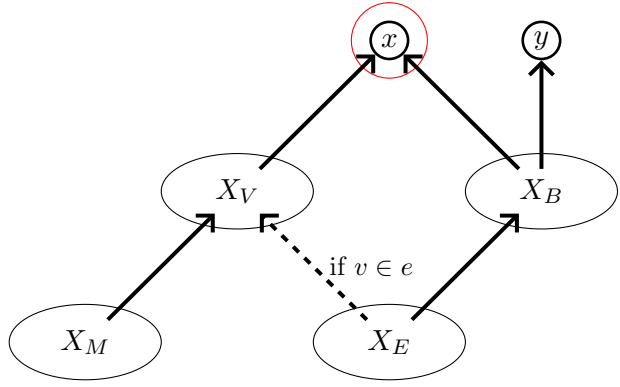
\begin{figure}
\begin{center}
\begin{tikzpicture}

 \draw (0,0) ellipse (1cm and 0.5cm) node(V) {$X_V$};
 \draw (-2,-2) ellipse (1cm and 0.5cm) node(M) {$X_M$};
 \draw (2,-2) ellipse (1cm and 0.5cm) node(E) {$X_E$};
 \draw (4,0) ellipse (1cm and 0.5cm) node(B) {$X_B$};

 \node[v](x) at (2,2) {$x$};
  \node[v](y) at (4,2) {$y$};
  \draw[draw=red, thin] (x) circle (0.5cm) ;

 \draw[-Straight Barb, line width=1.5pt] (M) -- (V);
 \draw[-Straight Barb, dashed, line width=1.5pt] (E) to node[right]{\small \,if $v\in e$} (V);
 \draw[-Straight Barb, line width=1.5pt] (V) -- (x);
 \draw[-Straight Barb, line width=1.5pt] (B) -- (x);
 \draw[-Straight Barb, line width=1.5pt] (B) -- (y);
 \draw[-Straight Barb, line width=1.5pt] (E) -- (B);

\end{tikzpicture}
\caption{\label{fig:redonesetsizeone} Reduction from {\sc Set Cover} to~\popoproblem. A solid arc means that all the possible relations are in the poset. The winning set is represented by a red circle.}
\end{center}
\end{figure}

\begin{proof}
  
  We reduce from the problem {\sc Set Cover}. A {\em set cover} of a hypergraph $H=(V,E)$ is a set $S \subseteq E$ of edges such that every vertex of $H$ belongs to some edge in $S$. Given an integer $k$, deciding if a given hypergraph contains a set cover of cardinality $k$ is {\sf NP}-complete~\cite{Karp1972}.

  Let $H=(V,E)$ be a hypergraph, and let $k$ be an integer. We assume that any vertex of $H$ is in at least one edge. We construct a \popogame~$\Ga_H=(X,\mathcal F, \poset)$ as follows. Define $m:=|E|$ and $n:=|V|$.
  The set of vertices $X$ is the union of:
  \begin{itemize}
  \item $X_V$, a set of size $n$, in one-to-one correspondence with $V$;
  \item $X_E$, a set of size $m$, in one-to-one correspondence with $E$;
  \item $X_M$, a set of size $m-k$;
  \item $X_B$, a set of size $n-k$;
  \item $x$ and $y$, two additional vertices.
  \end{itemize}

  There are $2(n+m-k+1)$ vertices in total.
  The set $\mathcal F$ contains a unique winning set $\{x\}$.
  As for the definition of $\poset$, we have $v<u$ in the following cases:
  \begin{itemize}
  \item $v\in X_M$ and $u\in X_V$;
  \item $v\in X_E$, $u\in X_V$, and $u\in v$ in $H$;
  \item $v\in X_E$, $u\in X_B$;
  \item $v\in X_V\cup X_B$, $u=x$;
  \item $v \in X_B$, $u=y$.
    \end{itemize}
 See Figure \ref{fig:redonesetsizeone} for an illustration of the structure of the game.
 Note that $\poset$ has height $3$. The predecessors of $x$ are in $X_V$ and $X_B$. Moreover, for all $v\in X_V$, $p(v)=2(n+m-k)$ is even, and for all $u\in X_B$, $p(u)=2(n+m-k)-1$ is odd. We will say a subset of $X$ is empty if all its vertices have been claimed: with Lemma~\ref{lem:predoneWS} in mind, Maker wants $X_V$ to be empty before $X_B$ is.

We will prove that $H$ has a set cover of size $k$ if and only if Breaker wins the game $\Ga_H$.

First of all, suppose that $H$ has a set cover of size $k$, $S=\{e_1,e_2,...,e_k\}$.
Both players can only claim vertices in $X_M \cup X_E$ until either $X_M$ or $X_E$ is empty.
During the first $m-k$ rounds, Breaker's strategy is the following. Whenever Maker claims a vertex in $X_E\setminus S$, Breaker claims an arbitrary vertex in $X_M$. Otherwise, Breaker claims an arbitrary vertex in $X_E\setminus S$.
This ensures that $X_E\setminus S$ is empty after exactly $m-k$ rounds.
Note that, among the first $m-k+1$ moves of Maker, she has necessarily claimed at least one vertex in $S$. Indeed, there are $2(m-k)$ vertices in $X_M\cup (X_E\setminus S)$, moreover any vertex in $X_V\cup X_B$ has a predecessor in $S$ and thus cannot be claimed while all vertices in $S$ are unclaimed. Therefore, after the $(m-k+1)^{\text{th}}$ move of Maker, $2(m-k)+1$ vertices have been claimed, all of which are in $X_M \cup X_E$, including all the vertices in $X_E \setminus S$ and at least one vertex in $S$.
This means that at most $k-1$ vertices in $X_E$ are unclaimed at this point. Breaker goes on to claim vertices in $X_E$ until it is empty (at most $k-1$ additional moves), and then moves on to claiming vertices in $X_B$ until is is empty (at most $n-k$ additional moves).
During this time, Maker cannot claim all $n$ vertices in $X_V$, so the last predecessor of $x$ to be claimed will be in $X_V$. By Lemma~\ref{lem:predoneWS}, Breaker wins.
\medskip

Now, suppose that $H$ has no set cover of size $k$. This means in particular that for any set $S \subseteq E$ of $i\leq k$ edges, at least $k-i+1$ vertices of $V$ are not covered by any edge in $S$. Indeed, one could otherwise add $k-i$ edges to cover all the vertices and obtain a set cover of size $k$.
We provide a strategy for Maker, in three stages. She starts by claiming an arbitrary vertex in $X_M$. She then answers the first $m-k-1$ moves of Breaker as follows:
\begin{itemize}
\item[(a1)] if Breaker claims a vertex in $X_E$, then Maker claims an arbitrary vertex in $X_M$;
\item[(a2)] otherwise, Breaker can only claim a vertex in $X_M$, and then Maker claims an arbitrary vertex in $X_E$.
  \end{itemize}

This concludes the first stage. At this point, $X_M$ is empty, and exactly $m-k-1$ vertices have been claimed in $X_E$. In the second stage, Maker answers the next $k+1$ moves of Breaker as follows:
\begin{itemize}
\item[(b1)] if Breaker claims a vertex in $X_E$, then Maker claims an arbitrary vertex in $X_V$;
\item[(b2)] otherwise, Breaker can only claim a vertex in $X_V$, and then Maker claims an arbitrary vertex in $X_E$.
\end{itemize}
First of all, note that Maker is always able to play the move in (b1). Indeed, consider that $i$ moves, $1\leq i\leq k$, have been played by Breaker in the second stage (including his last move in $X_E$). Then $k+1-i$ vertices remain unclaimed in $X_E$, and $n-i$  vertices are unclaimed in $X_V$. The $k+1-i$ unclaimed vertices in $X_E$ correspond to a set of edges of cardinality $k+1-i\leq k$. As we observed earlier, at least $i$ vertices of $V$ are not covered by this set and thus have all their predecessors claimed already. Of those vertices, at most $i-1$ have already been claimed, so at least one vertex in $X_V$ can be claimed by Maker.

After the second stage, $X_E$ is empty, $n-k-1$ vertices remain unclaimed in $X_V$, and $n-k$ vertices remain unclaimed in $X_B$. Therefore, Maker can claim all the remaining vertices in $X_V$ before Breaker claims all the vertices in $X_B$. By Lemma \ref{lem:predoneWS}, Maker wins.
\end{proof}

\begin{remark}
When there is only one winning set of size 1, the Maker-Maker convention is equivalent to the Maker-Breaker convention. Therefore, \popoproblem\ with only one winning set of size 1 and height 3 is also {\sf NP}-hard in Maker-Maker convention.
\end{remark}

\section{Posets of bounded width}\label{sec:width}
In this section, we consider posets of bounded width. For width 1, the game is completely determined and thus polynomial. However, for width 2, we prove that \popoproblem\ is already {\sf PSPACE}-hard even if the winning sets have size at most 3. However, if the number of winning sets is also bounded, we give an algorithm running in time $O(|X|^w2^mw^4)$.

\subsection{Width 2}

The reduction will be made from \QBF{}, which has been proven {\sf PSPACE}-complete by Stockmeyer and Meyer~\cite{stockmeyer1973}.\\

\noindent \QBF~\\
{\bf Input:}  A formula of the form $\psi = \forall x_1 \exists x_2 \dots \forall x_{2n-1} \exists x_{2n} \hspace{.2cm} \varphi(x_1,...,x_{2n-1})$ where $\varphi$ is a $3$-CNF formula \\
{\bf Output:} The truth value of $\psi$\\

\begin{theorem}\label{thm:psp_w2}
    \popoproblem\ is {\sf PSPACE}-complete even when restricted to instances where the poset is of width~2 and the winning sets are of size 3.
\end{theorem}

\begin{figure}
    \centering

\begin{tikzpicture}[rotate = 90]

\foreach \i in {1,...,4}
{
\draw ( 2 * \i , 0.5 ) node[v, inner sep=2](v\i){$v_{\i}$};
\draw ( 2 * \i , -1 ) node[v, inner sep=2](nv\i){$\overline{v_{\i}}$};
}
\foreach \i in {1,...,3}
{
\draw ( 2* \i + 1 , 2 ) node[v, inner sep=2](u\i){$u_{\i}$};
}

\draw[->] (v1) -- (u1);
\draw[->] (nv1) -- (u1);
\draw[->] (v2) -- (u2);
\draw[->] (nv2) -- (u2);
\draw[->] (v3) -- (u3);
\draw[->] (nv3) -- (u3);

\draw[->] (u1) -- (v2);
\draw[->] (u1) -- (nv2);
\draw[->] (u2) -- (v3);
\draw[->] (u2) -- (nv3);
\draw[->] (u3) -- (v4);
\draw[->] (u3) -- (nv4);

\draw[color= red, thin] (2,-0.15) arc (270:90:0.65) -- (6,1.15) arc (90:-90:0.65) -- (2,-0.15);
\draw[color= red, thin] (3.7,-0.6) -- (5.7,0.9) arc (126.9:53.1:0.5) -- (8.3,-0.6) arc (53.1:-126.9:0.5) -- (6.3, -0.35) arc (53.1:126.9:0.5) -- (4.3,-1.4) arc (306.9:126.9:0.5);
\end{tikzpicture}   
    \caption{The \popogame~obtained by reduction of $\forall x_1 \exists x_2 \forall x_3 \exists x_4 (x_1 \vee x_2 \vee x_3) \wedge (\neg x_2 \vee x_3 \vee \neg x_4)$}
    \label{fig:reduction width 2}
\end{figure}
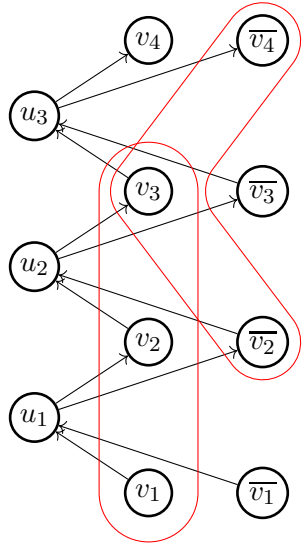

\begin{proof}

The result is obtained by reduction from \QBF. Let $\psi = \forall x_1 \exists x_2 \dots \forall x_{2n-1} \exists x_{2n} \varphi$ be a \QBF{} formula. 

A classical equivalent formulation of \QBF{} is to consider that two players, Satisfier and Falsifier, are choosing the valuations of the variables, with Falsifier starting and choosing the valuations of the odd variables whereas Satisfier chooses the valuations of the even variables. Satisfier wins if the valuation created by their moves satisfies $\varphi$, otherwise Falsifier wins.

We build a \popogame~$ \Ga = (X, \WS, \poset)$ as follows:

\begin{itemize}
    \item For any $1 \le i \le 2n$, we add two vertices $v_i$ and $\overline{v_i}$ in $X$.
    \item For any $1 \le i \le 2n-1$, we also add a vertex $u_i$.
    \item For any clause $C_j = l_{i_1} \vee l_{i_2} \vee l_{i_3}$ in $\varphi$, we add a winning set $S_j \in \WS$. For $1 \le k \le 3$, we have $v_{i_k} \in S_j$ if $l_{i_k} = x_{i_k}$ and $\overline{v_{i_k}} \in S_j$ if $l_{i_k} = \neg x_{i_k}$.
    \item For any $1 \le i \le 2n-1$, we add in $\poset$ the relations $v_i < u_i$, $\overline{v_i} < u_i$, $u_i < v_{i+1}$ and $u_i < \overline{v_{i+1}}$.
\end{itemize}

See Figure~\ref{fig:reduction width 2} for an example of the reduction. We now prove that Falsifier wins on $\psi$ if and only if Maker wins in $\Ga$.
First note that, because of the poset $\poset$, the moves are almost all forced. For any $1 \le i \le n$, Maker chooses one of $v_{2i-1}$ or $\overline{v_{2i-1}}$, then Breaker has to claim the other and Maker has to claim $u_{2i-1}$. Breaker then chooses one of $v_{2i}$ or $\overline{v_{2i}}$, Maker has to claim the other and Breaker has to claim $u_{2i}$, and so on.

Now suppose that Falsifier has a winning strategy $\strat$ on $\psi$. We define Maker's strategy as follows. For any $1 \le i \le n$, when Maker has to choose between $v_{2i-1}$ and $\overline{v_{2i-1}}$, she considers that all the variables $x_k$ with $k <  2i-1$ such that Maker has claimed $v_k$ are set to False, and the others to True. If, according to $\strat$, Falsifier has to put $x_{2i-1}$ to True, then Maker claims $\overline{v_{2i-1}}$, otherwise she claims $v_{2i-1}$. Since $\strat$ is a winning strategy for Falsifier, at the end of the game there exists a clause $C_j$ that is not satisfied by the valuation provided by Satisfier and Falsifier, i.e.~all its literals are false. By construction, the three corresponding vertices will be claimed by Maker, who will have claimed all the vertices of the winning set $S_j$.

Conversely, suppose Satisfier has a winning strategy $\strat$ on $\psi$. We define Breaker's strategy as follows. For any $1 \le i \le n$, when Breaker has to choose between $v_{2i}$ and $\overline{v_{2i}}$, he considers that all the variables $x_k$ with $k <  2i$ such that Maker has claimed $v_k$ are set to False, and the others to True. If, according to $\strat$, Satisfier has to put $x_{2i}$ to True, then Breaker claims $v_{2i}$, otherwise he claims $\overline{v_{2i}}$. Since $\strat$ is a winning strategy for Satisfier, at the end of the game every clause $C_j$ contains a literal that is true. By construction, the corresponding variable will be claimed by Breaker, who will have claimed at least one vertex in each $S_j$.
\end{proof}

\subsection{Parameterized algorithm with the number of winning sets}

In the following theorem, we give an algorithm of time complexity $O(|X|^w2^mw^4)$ for \popoproblem\ where $m$ is the number of winning sets and $w$ is the width of the poset.

\begin{theorem}\label{thm:xp_wl}
   \popoproblem\ can be decided in time $O(|X|^w2^mw^4)$ for instances where the poset is of width at most $w$ and there are at most $m$ winning sets.
\end{theorem}

\begin{proof}
We proceed by dynamic programming. Let $\Ga=(X,\WS,\poset)$ be a poset positional game.
Let $\WS=\{S_1, \ldots, S_m\}$ be the winning sets of $\Ga$.
We characterize the games that can be reached from a game played on $\Ga$ using an antichain $Y$ of $X$ and a boolean vector $B=(b_1,...,b_m)\in\{0,1\}^m$. We define the game $\Ga(Y,B)$ as follows:
\begin{itemize}
    \item The set of vertices of $\Ga(Y,B)$ is $X(Y,B)=\{x\in X | \,\exists y\in Y, y\leq x\}$.
    \item The set of winning sets of $\Ga(Y,B)$ is
    $\WS(Y,B)=\{S_i\cap X(Y,B) | \, 1 \leq i \leq m \text{ and } b_i=1\}$.
    \item The poset of $\Ga(Y,B)$ is $\poset$ restricted to $X(Y,B)$.
    \item As usual, Maker is the next player.
\end{itemize}

The set $Y$ corresponds to the available moves. The vector $B$ registers, for each winning set $S_i$, whether it contains an element that has already been claimed by Breaker (then $b_i=0$) or not (then $b_i=1$). In the latter case, it means that Maker has claimed all the elements of $S_i$ below $Y$, so she will win if she claims the other elements of $S_i$.
Note that the initial game $\Ga$ corresponds to $\Ga(Y_{\min},(1,...,1))$ where $Y_{\min}$ is the set of minimal elements of $\poset$. The outcome of $\Ga(Y,B)$ can be computed inductively as follows:
\begin{enumerate}
\item If $Y$ contains a winning set $S'_i\in \WS(Y,B)$ of size 1, then Maker wins. 
\item If all $b_i$'s are $0$, then Breaker wins.
\item Assume Conditions 1 and 2 are not met, and consider the next round. Let $u\in Y$ be a possible move for Maker. Note that Maker cannot win instantly with this move since condition 1 is not met. Let $Y^u$ be the set of minimal elements after $u$ is claimed, and let $v\in Y^u$ be a possible answer of Breaker. Finally, let $Y^{uv}$ be the set of minimal elements after $v$ is claimed. We define $B^{uv}\in \{0,1\}^m$ as follows: for each $1 \leq i \leq m$, let $b^{uv}_i=0$ if $v\in S_i$, otherwise let $b^{uv}_i=b_i$.
Then Maker wins $\Ga(Y,B)$ if and only if there exists $u\in Y$ such that, for all $v\in Y^u$, Maker wins $\Ga(Y^{uv},B^{uv})$.
\end{enumerate}

The algorithm consists in computing the outcome of all subgames, starting with $Y = \emptyset$, and traversing the subgames in increasing order of $|X(Y,B)|$.
Let us analyze the complexity of the algorithm.
There are at most $|X|^w$ antichains in a poset of width $w$ on $|X|$ elements, and at most $2^m$ possible values for $B$. Thus there at most $|X|^w2^m$ games $\Ga(Y,B)$ to consider.
Computing the outcome of
$\Ga(Y,B)$, when already knowing the outcome of
$\Ga(Y',B')$ for every antichain $Y'$ above $Y$, can be done in time $O(w^4)$. Indeed, there are $O(w^2)$ choices for the moves $u$ and $v$ of Maker and Breaker, and computing $Y^{uv}$ from $Y$ can be done in time $O(w^2)$. 
Consequently, the time complexity of the algorithm is
$O(|X|^w 2^mw^4)$.
\end{proof}

\begin{remark}
    Note the result can easily be adapted to the Maker-Maker convention. Indeed, it suffices to use elements $b_i$ with three values $\{0, 1, 2\}$: $b_i = 0$ (respectively 1) means that all elements of $S_i$ below $Y$ have been claimed by the first player (resp.~second player), and $b_i=2$ if both players have claimed vertices in $S_i$.
We obtain an algorithm in time $O(|X|^w 3^mw^4)$.
\end{remark}

\section{Posets made of pairwise disjoint chains}\label{sec:disjoint}

In this section, motivated by the game Connect-4, we consider posets made of pairwise disjoint chains. 
First note that if the size and the number of winning sets are bounded, then Theorem \ref{thm:xp_wl} applies.

\begin{corollary}\label{cor:disjointboundednumbersize}
\popoproblem on posets made of pairwise disjoint chains, with at mose $m$ winning sets and winning sets of size at most $s$ can be done in time $O(2^m|X|^{ms+1}(ms)^2)$.
\end{corollary}

\begin{proof}
    Since the number of winning sets is at most $m$ and there are at most $s$ elements in a winning set, there are at most $ms$ chains of $\poset$ that intersect a winning set. Let $Y$ be the set of all vertices which are in a chain of $\poset$ that does not intersect any winning set. In other words, the elements of $Y$ are not comparable with any element of a winning set.
    Let $X'=X\setminus Y$. Let $\Ga'=(X',\WS,\poset')$ where $\poset'$ is the poset induced by $\poset$ on $X'$.

     By Lemma \ref{lem:emptyUnion}, if $|Y|$ is even, then $\Ga$ and $\Ga'$ have the same outcome. If $|Y|$ is odd, then $\Ga$ has the same outcome than $\Ga''$, the game obtained from $\Ga'$ by adding a single element incomparable with the others.
    In both cases, we reduce the problem of deciding the outcome of $\Ga$ to the problem of deciding the outcome of a game with: $m$ winning sets, at most $|X|$ vertices, and a poset of width at most $ms+1$. By Theorem \ref{thm:xp_wl}, it can be solved in time $O(2^m|X|^{ms+1}(ms)^2)$.
\end{proof}

In the rest of this section, we give positive results for \popogames~on pairwise disjoint chains with very small winning sets (but an unbounded number of them). We will often use parity arguments when discussing games on posets made of pairwise disjoint chains. Therefore, we introduce the following framework, which will be useful in this section. Consider a \popogame~$\Ga=(X,\WS,\poset)$ where $\poset$ is made of pairwise disjoint chains $C_1,...,C_w$. The elements of the chain $C_i$ are denoted as $x_{i,1}>x_{i,2}>...>x_{i,\ell_i}$. This numbering of the vertices from top to bottom is best adapted to define the following coloring: a vertex $x_{i,j}$ is colored white if $j$ has same parity as $|X|$, otherwise it is colored black. Note that this coloring satisfies the following properties:
\begin{itemize}
    \item The coloring does not change after a full round. When Breaker is playing, we keep the same coloring. 
    \item All the maximal elements are of the same color, that is: white if $|X|$ is odd, or black if $|X|$ is even.
    \item If $|X|$ is even, then the number of black vertices that are minimal is even (resp.~odd) when it is Maker's (resp.~Breaker's) turn.
    \item If $|X|$ is odd, then the number of white vertices that are minimal is odd (resp.~even) when it is Maker's (resp.~Breaker's) turn.
\end{itemize}

We first establish some elementary lemmas about winning sets of size 1, which help provide a linear-time algorithm in the case where all winning sets are of size 1. Recall that, by Remark \ref{rem:ws1k}, this is equivalent to considering a single winning set of unbounded size: we use this to solve some instances of \textit{Connect-$k$} in the Maker-Breaker convention. 
Secondly, we use dynamic programming to design an algorithm running in time $O(h^{2(w+1)})$ when all the winning sets have size at most 2 (where $h$ and $w$ respectively denote the height and the width of the poset).
Finally, we turn to the case where all winning sets are of size at most 2 and all chains are of height at most 2, for which we give an algorithm in time $O(|X|^4)$.

\subsection{Winning sets of size 1}

The first observation is that Maker wins if there is a winning set of size 1 which is white.

\begin{lemma}\label{lem:white}
    If $\{u\}\in \WS$ for some white vertex $u$, then Maker wins.
\end{lemma}

\begin{proof}
Write $u=x_{i,j}$. If $x_{i,j}$ is a minimal element of $\poset$, then Maker can claim it as her first move. Otherwise, consider its unique predecessor $x_{i,j+1}$.
The number of vertices that are greater or equal than $x_{i,j}$ is exactly $j$, which has same parity as $|X|$ since $x_{i,j}$ is colored white.
Thus, the number of vertices that are not greater or equal than $x_{i,j}$, i.e.~that can be claimed before $x_{i,j}$, is even. Therefore, if Maker never claims $x_{i,j+1}$, then Breaker will have to claim it, after which Maker can claim $x_{i,j}$ and win.
\end{proof}

It is easy to see that the same assertion would be false for black vertices. However, if $|X|$ is odd and there are two such winning sets on different chains, then Maker is also winning.

\begin{lemma}\label{lem:blackblack}
    If $|X|$ is odd and $\{u\},\{v\} \in \WS$ for some black vertices $u$ and $v$ sitting on different chains, then Maker wins.
\end{lemma}

\begin{proof}
If $u$ or $v$ is a minimal element, then the result is obvious. Otherwise, let $u'$ (resp.~$v'$) be the predecessor of $u$ (resp.~of $v$), and let $Y_1$ (resp.~$Y_2$) be the set of all vertices that are greater or equal than $u'$ (resp.~than $v'$). Since $|X|$ is odd and $u',v'$ are white, we know $|Y_1|$ and $|Y_2|$ are odd, so $|X \setminus (Y_1 \cup Y_2)|$ is odd. Maker's strategy is as follows. She starts by claiming an arbitrary vertex in $X \setminus (Y_1 \cup Y_2)$. Whenever Breaker claims a vertex in $X \setminus (Y_1 \cup Y_2)$, Maker answers by claiming another arbitrary vertex in $X \setminus (Y_1 \cup Y_2)$. The parity ensures that Breaker will have to claim either $u'$ or $v'$, at which point Maker claims its successor and wins.
\end{proof}

It turns out these two elementary lemmas are sufficient to solve the general case where all the winning sets are of size 1.

\begin{theorem}\label{thm:disjoint-size1}
    Let $\Ga=(X,\WS,\poset)$ be a \popogame~where $\poset$ is made of pairwise disjoint chains and all elements of $\WS$ are of size 1. Maker wins $\Ga$ if and only if at least one of the following conditions hold:
    \begin{itemize}
        \item there is a winning set that is a minimal element; or
        \item there is a white winning set; or
        \item $|X|$ is odd and there are two black winning sets on different chains.
    \end{itemize}
    In particular, \popoproblem\ can be solved in linear time for such games.
\end{theorem}

\begin{proof}
    The ``if'' direction is a direct consequence of Lemmas~\ref{lem:white}~and~\ref{lem:blackblack}. Let us prove the ``only if'' direction: we suppose that none of the three conditions hold, and we show that Breaker wins. In particular, all winning sets are black (we may assume that there exists at least one) and no minimal element is a winning set. Additionally, we are in one of two cases:
    \begin{itemize}
	   \item Case 1: $|X|$ is even. This means there is an even number of chains of odd size, i.e.~$|Y|$ is even where $Y$ denotes the set of all black minimal elements. Breaker's strategy is the following. Whenever Maker claims a vertex in $Y$, Breaker also claims an arbitrary vertex in $Y$. Otherwise, Maker necessarily claims a white vertex, so Breaker claims its successor (which always exists since all maximal elements are black). Using this strategy, Breaker ensures to claim every black vertex in $X \setminus Y$, so he wins.
	   \item Case 2: $|X|$ is odd and all the winning sets are part of the same chain $C_i$. Let $j$ be maximum such that $\{x_{i,j}\}\in\WS$, and let $Y$ be the set of all vertices which are not greater or equal than $x_{i,j+1}$. Since $|X|$ is odd and $x_{i,j+1}$ is white, we know $|Y|$ is even. Breaker's strategy is the following. Whenever Maker claims a vertex in $Y$, Breaker also claims an arbitrary vertex in $Y$. Otherwise, Maker necessarily claims a white vertex in $C_i$, so Breaker claims its successor (unless Maker has claimed $x_{i,1}$ specifically, in which case Breaker plays an arbitrary move). Using this strategy, Breaker ensures to claim every black vertex in $X \setminus Y$, so he wins. \qedhere
    \end{itemize}
\end{proof}

Using this result, we can solve a basic case of \textit{$w \times h$ Connect-$k$}  in Maker-Breaker convention, which is the \popogame~played on $w$ chains of height $h$ where the winning sets are the alignments of $k$ consecutive vertices (horizontal, vertical or diagonal). Note that the game is easily solved if $k \leq 2$, or if $k>w$, or if $k>h$ and $wh$ is even.

\begin{corollary}\label{cor:connect-k}
    Let $k,w,h$ be integers with $3 \leq k \leq w$ and $h > 1$. If $w$ and $h$ are both odd, then Maker wins the $w \times h$ Connect-$k$ game in Maker-Breaker convention.
\end{corollary}

\begin{proof}
    We actually show that a single well-chosen winning set if enough for Maker to win. Let $W$ be a horizontal alignment of $k$ vertices on the top row: we claim that Maker wins $\Ga=(X,\WS,\poset)$ where $X$ and $\poset$ are defined as in $w \times h$ Connect-$k$ and $\WS=\{W\}$. For this, using the idea from Remark~\ref{rem:ws1k}, she claims an arbitrary vertex $x$ and then sees herself as Breaker (playing second, as usual) in the game $\Ga'=(X',\WS',\poset')$ where: $X'=X \setminus \{x\}$, $\poset'$ is the poset induced by $\poset$ on $X'$, and $\WS'=\{\{y\},y \in W\}$. Since $h \geq 3$, all vertices of $W$ are nonminimal in $\poset'$. Moreover, they are maximal in $\poset'$, so they are black since $|X'|=wh-1$ is even. Therefore, Theorem~\ref{thm:disjoint-size1} ensures that Breaker wins $\Ga'$, so Maker wins $\Ga$.
\end{proof}

\begin{remark}
    It should be noted that Connect-$k$, unlike all other game classes studied in this paper, requires a separate study of the case where Breaker plays first, since this case does not reduce to the case where Maker plays first (after Breaker's first move, we do not have an instance of Connect-$k$ anymore). However, from the same argument used in the proof of Corollary~\ref{cor:connect-k}, we get that Maker wins playing second if $wh$ is even. 
\end{remark}

\subsection{Winning sets of size at most 2}

Things become more complicated when the winning sets are not all of size 1. We are going to need the following reduction, using a similar idea to the proof of Theorem \ref{thm:disjoint-size1}. The point is that, in the presence of a black winning set of size 1, we can simplify the game (recall that Lemma~\ref{lem:white} already solves the case where there is a white winning set of size 1). Using the same notations as before, assume that there exists a winning set $\{x_{i,j}\}$ for some nonminimal black vertex $x_{i,j}$ with $j \geq 3$.
We define a reduced game $\Ga'=(X',\WS',\poset')$ as follows.

\begin{figure}[ht]
\begin{center}
\scalebox{0.7}{    \begin{tikzpicture}
        
    \node[v, inner sep=2](11) at (0,10.5) {$x_{1,1}$};
    \node[v, fill=black!50, inner sep=2](12) at (0,9) {{\color{white} $x_{1,2}$}};
    \node[v, inner sep=2](13) at (0,7.5) {$x_{1,3}$};
    \node[v, fill=black!50, inner sep=2](14) at (0,6) {{\color{white} $x_{1,4}$}};
    \node[v, inner sep=2](15) at (0,4.5) {$x_{1,5}$};
    \node[v, fill=black!50, inner sep=2](16) at (0,3) {{\color{white} $x_{1,6}$}};
    \node[v, inner sep=2](17) at (0,1.5) {$x_{1,7}$};

    \draw[->] (17) -- (16);
    \draw[->] (16) -- (15);  
    \draw[->] (15) -- (14);  
    \draw[->] (14) -- (13);  
    \draw[->] (13) -- (12);  
    \draw[->] (12) -- (11);

    \node[v, inner sep=2](21) at (3,10.5) {$x_{2,1}$};
    \node[v, fill=black!50, inner sep=2](22) at (3,9) {{\color{white} $x_{2,2}$}}; 
    \node[v, inner sep=2](23) at (3,7.5) {$x_{2,3}$};
    \node[v, fill=black!50, inner sep=2](24) at (3,6) {{\color{white} $x_{2,4}$}};
    \node[v, inner sep=2](25) at (3,4.5) {$x_{2,5}$};
    \node[v, fill=black!50, inner sep=2](26) at (3,3) {{\color{white} $x_{2,6}$}};     
  
    \draw[->] (26) -- (25);  
    \draw[->] (25) -- (24);  
    \draw[->] (24) -- (23);  
    \draw[->] (23) -- (22);  
    \draw[->] (22) -- (21); 

    \draw[color= red, thin] (16) circle (0.65cm);
    \draw[color= red, thin] (0,8.15) arc (90:270:0.65) -- (3,6.85) arc (-90:90:0.65) -- (0,8.15);

  \draw[color= red, thin] (0,5.35) arc (270:90:0.65) arc (270:326.3:0.65) -- (1.92,9.07) arc (146.3:90:1.3);
    \draw[color= red, thin] (3,9.65) arc (90:-90:0.65) arc (90:146.3:0.65) -- (1.08,5.93) arc (326.3:270:1.3);

     \draw[color= red, thin] (0,0.85) arc (270:90:0.65) arc (270:333.4:0.7) -- (1.75,5.87) arc (153.4:90:1.4);
       \draw[color= red, thin] (3,6.65) arc (90:-90:0.65) arc (90:153.4:0.7) -- (1.25,1.63) arc (333.4:270:1.4);

\draw[color=blue,dashed, rounded corners] (-0.8,5.2) rectangle (0.7,11.3);
\node at (-1.2,8) {\huge \color{blue} $Y$};

\node at (-2.5,6) {\huge $\Ga$};

\node at (6,6) {\Huge $\Leftrightarrow$};

\begin{scope}[shift={(9,0)}]
    \node[v, inner sep=2](15) at (0,4.5) {$x_{1,5}$};
    \node[v, fill=black!50, inner sep=2](16) at (0,3) {{\color{white} $x_{1,6}$}};
    \node[v, inner sep=2](17) at (0,1.5) {$x_{1,7}$};

    \draw[->] (17) -- (16);
    \draw[->] (16) -- (15);

    \node[v, inner sep=2](21) at (3,10.5) {$x_{2,1}$};
    \node[v, fill=black!50, inner sep=2](22) at (3,9) {{\color{white} $x_{2,2}$}}; 
    \node[v, inner sep=2](23) at (3,7.5) {$x_{2,3}$};
    \node[v, fill=black!50, inner sep=2](24) at (3,6) {{\color{white} $x_{2,4}$}};
    \node[v, inner sep=2](25) at (3,4.5) {$x_{2,5}$};
    \node[v, fill=black!50, inner sep=2](26) at (3,3) {{\color{white} $x_{2,6}$}};     
  
    \draw[->] (26) -- (25);  
    \draw[->] (25) -- (24);  
    \draw[->] (24) -- (23);  
    \draw[->] (23) -- (22);  
    \draw[->] (22) -- (21); 

    \draw[color= red, thin] (16) circle (0.65cm);
    \draw[color= red, thin] (23) circle (0.65cm);

    \draw[color= red, thin] (0,0.85) arc (270:90:0.65) arc (270:333.4:0.7) -- (1.75,5.87) arc (153.4:90:1.4);
    \draw[color= red, thin] (3,6.65) arc (90:-90:0.65) arc (90:153.4:0.7) -- (1.25,1.63) arc (333.4:270:1.4);
    
    \node at (5.5,6) {\huge $\Ga'$};

\end{scope}
\end{tikzpicture}}
\end{center}
\caption{\label{fig:reduc} Reduction of a game $\Ga$ containing a black winning set of size 1 to a game $\Ga'$. The winning set $\{x_{1,4},x_{2,2}\}$ disappears since it contains a black vertex in $Y$. The winning set $\{x_{1,3},x_{2,3}\}$ becomes the winning set $\{x_{2,3}\}$ in $\Ga'$ since $x_{1,3}$ is white. By Lemma~\ref{lem:reducWS1}, the two games have the same outcome. In this case, since $\Ga'$ has a white winning set of size 1, both games are winning for Maker.}
\end{figure}
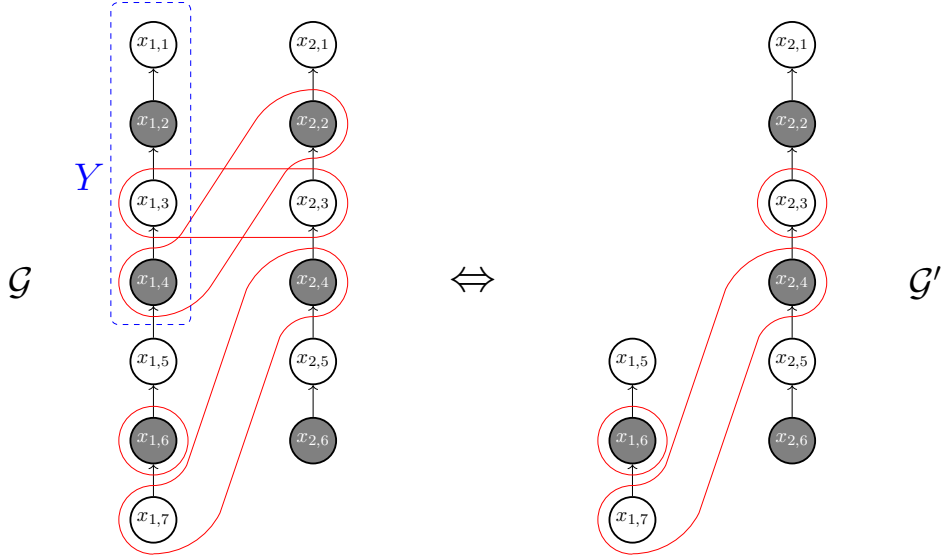

If $|X|$ is even, then define $Y=\{x_{i,1},...,x_{i,j-1}\}$, otherwise define $Y=\{x_{i,1},...,x_{i,j-2}\}$. In both cases, let $X'=X\setminus Y$ and let $\poset'$ be the poset induced by $\poset$ on $X'$. Note that, since $x_{i,j}$ is black, we have removed an even number of vertices in both cases. In particular, the coloring of the vertices is the same for $(X',\poset')$ as for $(X,\poset)$.
 We now define the winning sets of $\Ga'$ as follows. Let $S\in \WS$. If $S\subseteq X'$, then $S\in \WS'$. If $S$ only intersects $Y$ on white vertices, then $S\cap X' \in \WS'$. Otherwise, i.e.~if $S$ contains a black vertex greater than $x_{i,j}$, then we ignore $S$. See Figure~\ref{fig:reduc} for an illustration in the case $|X|$ odd. We will prove that $\Ga$ and $\Ga'$ have the same outcome. For this, we first need some technical lemmas.

\begin{lemma}\label{lem:ordervertices}
If $|X|$ is even and Maker has a winning strategy in $\Ga'$, then she has a winning strategy in $\Ga'$ such that she never claims $x_{i,j+1}$ unless it makes her win on the spot.
\end{lemma}

\begin{proof}
Let $\strat$ be a winning strategy for Maker in $\Ga'$. If, at any point during the game, Maker claims $x_{i,j+1}$ without winning instantly, then Breaker is forced to claim $x_{i,j}$, otherwise Maker would win on her next turn.
Let $\strat'$ be the same strategy as $\strat$ with the exception that, when Maker should claim $x_{i,j+1}$, she skips this claim (and the forced answer $x_{i,j}$ of Breaker) and goes on with the game, following the strategy $\strat$ as if that round had happened. Breaker will never claim $x_{i,j+1}$, otherwise he would lose immediately. Note that, since $x_{i,j+1}$ is white, there will always be another vertex for Breaker to claim.
At some point in $\strat$, Maker will fill up a winning set $W \in \WS'$. If $x_{i,j+1}\notin W$, then $W$ will also be filled up in $\strat'$. If $x_{i,j+1}\in W$, then at some point in $\strat'$ Maker will have claimed all the vertices of $W$ except for $x_{i,j+1}$: in this case, Maker wins on the spot by claiming $x_{i,j+1}$, satisfying the statement of the lemma.
\end{proof}

We give a similar result for Breaker.

\begin{lemma}\label{lem:orderverticesbreaker}
    If $|X|$ is odd and Breaker has a winning strategy in $\Ga'$, then he has a winning strategy in $\Ga'$ where he never claims $x_{i,j-1}$.
\end{lemma}

\begin{proof}
Consider a winning strategy $\strat$ for Breaker in $\Ga'$. Note that Breaker never claims $x_{i,j+1}$ in $\strat$, otherwise he would lose immediately. We now define a new strategy $\strat'$ as follows. Breaker plays according to $\strat$ until, inevitably, Maker claims $x_{i,j+1}$: at this point, Breaker claims $x_{i,j}$, but he then goes on with $\strat$ as if that last round had not happened. In particular, Breaker will never be told to claim $x_{i,j-1}$, as this vertex is not even available in this narrative. Therefore, Maker will eventually claim $x_{i,j-1}$: only at this moment does Breaker input into $\strat$ Maker's claims of $x_{i,j+1}$ (prompting the answer $x_{i,j}$ of Breaker) and $x_{i,j-1}$ successively. Breaker then follows $\strat$ until the end, winning the game without claiming $x_{i,j-1}$.
\end{proof}

\begin{lemma}\label{lem:reducWS1}
Assume that $\Ga$ has no winning set containing only white vertices that are greater than $x_{i,j}$.
Then the two games $\Ga$ and $\Ga'$ have the same outcome when Maker starts.
\end{lemma}

\begin{proof}



Firstly, suppose Maker has a winning strategy $\strat'$ in $\Ga'$. If $|X|$ is even, then we further assume that $\strat'$ satisfies the property from Lemma \ref{lem:ordervertices}. Up to playing arbitrary moves after filling up a winning set, assume that $\strat'$ continues until all vertices in $X'$ are claimed. The idea is the following: Maker will fill up a winning set of $\Ga'$, and also claim all the white vertices in $Y$. By definition of $\Ga'$, this will automatically ensure that Maker fills up a winning set of $\Ga$. We now explain how she achieves this:
\begin{itemize}
    \item Case 1: $|X|$ is odd. Maker plays her first move $x \in X'$ according to $\strat'$. Since $|X' \setminus \{x\}|$ and $|Y|$ are both even, Maker can then use the following strategy. Whenever Breaker plays inside $X'$, Maker answers inside $X'$ according to $\strat'$, thus ensuring to eventually fill a winning set of $\Ga'$.
    Otherwise i.e.~whenever Breaker claims some (black) $y \in Y$, Maker claims the (white) successor of $y$, thus ensuring to claim all the white vertices in $Y$.
    \item Case 2: $|X|$ is even. Maker plays according to $\strat'$ until $x_{i,j+1}$ is claimed. Assume it is claimed by Maker, otherwise she wins immediately by claiming $x_{i,j}$. By the property from Lemma \ref{lem:ordervertices}, Maker has already filled up a winning set of $\Ga'$ at this point, so she just needs to claim all the white vertices of $Y$. Breaker is forced to claim $x_{i,j}$, and then Maker claims $x_{i,j-1}$. From this point on, since $j-2$ is odd, Maker can force Breaker to claim all the black vertices in $Y$, and claim all their (white) successors herself.
\end{itemize}




Lastly, suppose Breaker has a winning strategy $\strat'$ in $\Ga'$. If $|X|$ is odd, then we further assume that $\strat'$ satisfies the property from Lemma \ref{lem:orderverticesbreaker}. The idea is the following: Breaker will claim a vertex in each winning set of $\Ga'$, and also claim all the black vertices in $Y$. By definition of $\Ga'$, this will automatically ensure that Breaker claims a vertex in each winning set of $\Ga$. We now explain how he achieves this:
\begin{itemize}
    \item Case 1: $|X|$ is even. Since $|X'|$ and $|Y|$ are both even, Breaker can then use the following strategy. Whenever Maker plays inside $X'$, Breaker answers inside $X'$ according to $\strat'$, thus ensuring to claim a vertex in each winning set of $\Ga'$. Otherwise i.e.~whenever Maker claims some (white) $y \in Y$, Breaker claims the (black) successor of $y$, thus ensuring to claim all the black vertices in $Y$.
    \item Case 2: $|X|$ is odd. Breaker plays according to $\strat'$ until Maker claims $x_{i,j-1}$, which will necessarily happen by the property from Lemma \ref{lem:orderverticesbreaker}. Breaker pretends like that move did not happen i.e.~does not input that claim into $\strat'$, and claims $x_{i,j-2}$ himself. At this point, the number of unclaimed vertices in $X'$ is even, so that, whenever Maker plays inside $X'$, Breaker can answer inside $X'$ according to $\strat'$. This will force Maker to claim all the white vertices in $Y \setminus\{x_{i,1}\}$, with Breaker claiming the (black) successor each time. Eventually, Maker will have to claim $x_{i,1}$, at which point Breaker inputs Maker's earlier claim of $x_{i,j-1}$ into $\strat'$. The sets of vertices claimed by both players in $\strat'$ now coincide with reality, which ensures that Breaker claims a vertex in each winning set of $\Ga'$. \qedhere
\end{itemize}
\end{proof}

Now equipped with this reduction, we turn to the case where all winning sets have size at most 2.

\begin{theorem}\label{thm:disjoint-size2}
    \popoproblem\ can be solved in time $O(h^{2(w+1)})$ for instances where the poset consists of $w$ pairwise disjoint chains of height at most $h$ and all winning sets have size at most 2.
\end{theorem}

\begin{proof}
Our algorithm is based on a dynamic programming approach. 
As before, we consider a \popogame\ $\Ga=(X,\WS,\poset)$ where $\poset$ is made of pairwise disjoint chains. We use the same notations and colouring than before.
Thanks to Lemmas \ref{lem:white} and \ref{lem:reducWS1}, we can reduce the number of sub-positions to consider. Indeed, we can assume that each chain contains at most one (black) winning set of size 1. Our algorithm computes dynamically the outcome of all these ``useful" positions.
A sub-position will be given by two $w$-uplets : $K=(k_1,...,k_w)$ with $0\leq k_i\leq \ell_i$ corresponding to the heights of the chains in the sub-position, and $W=(s_1,...,s_w)$ with $0\leq s_i \leq k_i$ corresponding to the index of a black winning set of size 1 on the $i$-th chain that might be created during the game. We will only consider values of $s_i$ that either are zero (no winning set of size 1 on the $i$-th chain) or do not have the same parity as $|X|$ (so that $x_{i,s_i}$ is necessarily black) and such that there is a winning set $\{x_{i,s_i},u\}$ in $\WS$.
The game $\Ga(K,W)$ is formally defined as follows:
\begin{itemize}
    \item The vertices $X(K,W)$ of $\Ga(K,W)$ on the $i$-th chain are:
    \begin{itemize}
        \item if $s_i=0$: the vertices $x_{i,1}$ to $x_{i,k_i}$,
        \item if $s_i>0$ and $|X|$ is even: the vertices $x_{i,1}$ to $x_{i,s_i-1}$,
        \item if $s_i>0$ and $|X|$ is odd: the vertices $x_{i,1}$ to $x_{i,s_i-2}$. 
    \end{itemize}
    \item The poset $\poset(K,W)$ is the poset induced by $\poset$ on $X(K,W)$.
    \item The set $\WS(K,W)$ contains all sets $\{x_{i,s_i}\}$ where $s_i\neq 0$, and additionally, for each $S\in \WS$:
    \begin{itemize}
        \item If $S$ has an element $x_{i,j}$ with $j>k_i$, then $S$ is ignored.
        \item If $S$ has a black element $x_{i,j}$ with $j<s_i$, then $S$ is ignored.
        \item Otherwise, if $|X|$ is even (resp. odd), we add $S \setminus \{x_{i,j}, j<s_i\}$ (resp. $S \setminus \{x_{i,j}, j<s_i-1\}$) to the set of winning sets $\WS(K,W)$. By doing this, we might add an empty winning set to $\WS(K,W)$, which actually means that Maker wins trivially.
    \end{itemize}
    \item As usual, Maker is the next player.
\end{itemize}

Note that initial game $\Ga$ corresponds to the game $\Ga(K,W)$ with $K=(\ell_1,...,\ell_w)$ and $W=(0,...,0)$. 
We compute the outcome in increasing order of the total number of vertices.

Let $\Ga(K,W)$ be a position. We first list some trivial cases, using Lemmas~\ref{lem:white}~and~\ref{lem:blackblack}:
\begin{itemize}
    \item If there are no winning sets, then Breaker wins.
        \item If there is an empty winning set, then Maker wins.
    \item If there is a winning set of size 1 which is white or minimal, then Maker wins.
    \item If there is a winning set made of two white vertices on the same chain, then Maker wins.
    \item If there are two black winning sets of size 1 on different chains and $|X|$ is odd, then Maker wins.
\end{itemize}
In particular, that last item applies as soon as there are two positive $s_i$'s and $|X|$ is odd. The fourth item comes from the fact that Maker can apply Lemma \ref{lem:white} twice in a row to claim both vertices of the winning set.

Assuming none of these trivial cases hold, we consider every possible move $x_{i,j}$ of Maker and every possible answer $x_{i',j'}$ of Breaker in the first round. Each choice leads to a sub-position $\Ga(K',W')$. We define $K'=(k'_1,\ldots,k'_w)$ as follows: $k'_i=k_i-1$, $k'_{i'}=k_{i'}-1$ (or, if $i=i'$: $k'_i=k_i-2$), and all the others values in $K'$ are equal to the corresponding values in $K$. As for $W'=(s'_1,...,s'_w)$, its definition depends on $x_{i,j}$ and $x_{i',j'}$.
\begin{itemize}
    \item Suppose $x_{i,j}$ is in no winning set.  Then $W'=W$ except that we define $s'_{i'}=0$ if $j'=s_{i'}$.
    \item Suppose $x_{i,j}$ is in at least one winning set. Let $U$ be the set of all vertices $u$ such that $\{x_{i,j},u\}\in\WS(K,W)$ and $u$ is not claimed by Breaker (i.e.~$u\neq x_{i',j'}$). Then, for each $u\in U$, $\{u\}$ is a winning set after the first round. We assume that $S$ only contains black vertices, otherwise Maker wins. For each $i''$, let $s'_{i''}$ be the maximum between $s_{i''}$ and all the $j''$ such that $x_{i'',j''} \in U$.
    Lemma \ref{lem:reducWS1} ensures that, with this definition of $W'$, the game $\Ga(K',W')$ is indeed equivalent to this case.
\end{itemize}

Knowing all the possible outcomes after these two moves, one can compute the outcome of $\Ga(K,W)$. If there exists a move of Maker such that, for any Breaker answer, the resulting game $\Ga(K',W')$ is winning for Maker, then Maker wins. Otherwise, Breaker wins.

There are $(h+1)^{w}$ possible values for $K$ and also $(\lceil h/2 \rceil+1)^{w}$ possible values for $W$ (since the winning sets of size 1 must be black) hence $O(h^{2w}/2^{w})$ positions to consider. To compute the outcome of a position, one needs to consider $O(w^2)$ potential pairs of moves. Each given pair may require to go through the whole set of winning sets, which is of size $|\WS|=O((hw)^2)$. All in all, the running time of the algorithm is $O(h^{2w}/2^w*w^2*h^2w^2)=O(h^{2(w+1)})$.
\end{proof}

\begin{remark}
    Thanks to Lemma \ref{lem:blackblack}, it is possible to reduce the number of sub-positions when $|X|$ is odd. Indeed, we can assume that there is at most one winning set of size $1$ in this case. Therefore, $W$ could be replaced by an ordered pair $(i,j)$ such that $\{x_{i,j}\}$ is the unique winning set of size 1, with the convention that $i=0$ if there is no such winning set. Then there will only be at most $wh/2$ possible values for $W$, hence a final complexity in $O(h^{w+3}w^5)$.
\end{remark}

\subsection{Winning sets of size at most 2, chains of height at most 2}

The algorithm from Theorem~\ref{thm:disjoint-size2} is not efficient for posets of unbounded width. However, when restricting the problem to posets of height at most 2, we do get a polynomial-time algorithm in this case also. Since the chains are of height at most 2, we will simply use the ``top/bottom" terminology already adopted in the proof of Theorem~\ref{thm:topverticesWS}, rather than the coloring previously used in this section.


\begin{theorem}\label{thm:poly_disj_d2}
    \popoproblem\ can be solved in time $O(|X|^4)$ for instances where the poset on $X$ consists of pairwise disjoint chains of height at most 2 and all winning sets are of size at most 2.
\end{theorem}

\begin{proof}
   
Let $x_1,\ldots,x_t,y_1,\ldots,y_t$ be the vertices of the chains of height 2, with $x_i<y_i$ for all $1 \leq i \leq t$. Let $x_{t+1},\ldots,x_w$ be the other vertices if there are any. We say the $x_i$ are the ``bottom" vertices and the $y_j$ are the ``top" vertices. We consider two cases depending on the parity of $|X|$. 

\begin{itemize}
    \item Firstly, suppose $|X|$ is odd, which means Maker plays the first and the last move. There are several easy winning situations for Maker (note that the first three actually hold regardless of parity), all of which can clearly be identified in time $O(|X|^2)$ since there are $O(|X|^2)$ winning sets: 
        \begin{enumerate}
            \item There is a winning set of the form $\{x_i\}$. Indeed, Maker simply claims $x_i$ as her first move.
            \item There are two winning sets of the form $\{x_i,x_j\}$ and $\{x_j,x_k\}$. Indeed, Maker starts by claiming $x_j$ and will claim either $x_i$ or $x_k$ as her second move to win.
            \item There are two winning sets of the form $\{x_i,x_j\}$ and $\{x_i,y_i\}$. Indeed, Maker starts by claiming $x_i$ and will claim either $x_j$ or $y_i$ as her second move to win.
            \item There is a winning set of the form $\{x_i,y_j\}$ with $i\neq j$. Indeed, Maker starts by claiming $x_i$, the parity then allows her to wait until Breaker claims $x_j$, at which point she claims $y_j$ and wins. 
            \item There is a winning set in which each vertex is a top vertex. Indeed, using the parity, Maker can ensure to claim all the top vertices in the game. For this, she starts by claiming $x_{t+1}$ (that exists since the total number of vertices is odd). Each time Breaker claims some $x_i$, Maker claims $y_i$ if $i \leq t$ or some $x_j$ with $j>t$ if $i>t$. 
        \end{enumerate}
        We actually claim that those are the only winning situations for Maker. Indeed, suppose none of the above apply: then the winning sets are pairwise disjoint and all of them are of the form $\{x_i,x_j\}$ or $\{x_i,y_i\}$, so whenever Maker claims a vertex of a winning set it is possible for Breaker to immediately claim the other.

    \item Secondly, suppose $|X|$ is even, which means Breaker plays the last move. We then introduce the following reduction rules:
    \begin{itemize}
        \item[R1:] Assume there exists a winning set of the form $\{x_i,x_j\}$ that is disjoint from all the other winning sets. Then remove $x_i$ and $x_j$.
        \item[R2:] Assume there exists a nonempty set of indices $I \subseteq \{1,\ldots,t\}$ such that all the winning sets that contain some $x_i$ with $i\in I$ also contain some $y_j$ with $j\in I$. Then remove all vertices with indices in $I$.
    \end{itemize}
    \
    
    $\quad$ Claim 1: The reductions R1 and R2 are outcome-neutral when Maker is the first player.
    
    \
    Let us prove Claim 1 (note that the proof heavily relies on the parity assumption):
    \begin{itemize}
        \item[--] Let $\Ga'$ be a reduced game obtained from $\Ga$ through R1. If Maker wins $\Ga'$, then she wins $\Ga$ by claiming $x_i$ as her first move: indeed, this forces Breaker to claim $x_j$, and we get the game $\Ga'$. If Breaker wins $\Ga'$, then he wins $\Ga$ by following his winning strategy in $\Ga'$ whenever Maker does not play inside $\{x_i,x_j\}$, and only claiming $x_i$ or $x_j$ when Maker claims the other. However, we need to be careful that Breaker might be forced to claim $x_i$ or $x_j$ before Maker does. Because of the parity assumption, this can only happen if, say, $i \leq t$ and $j>t$, and $x_i,y_i,x_j$ are the only three vertices remaining. In this case, Breaker simply claims $x_j$: Maker must claim $x_i$ and Breaker claims $y_i$, ensuring the win since $x_i$ is in no other winning set.
        \item[--] Let $\Ga'$ be a reduced game obtained from $\Ga$ through R2, and let $X_I$ be the set of vertices with indices in $I$. Then whoever wins $\Ga'$ also wins $\Ga$, by following their winning strategy in $\Ga'$ whenever their opponent does not play inside $X_I$, and claiming $y_i$ whenever their opponent claims $x_i$ for some $i \in I$. This strategy clearly works for Maker, who is going to fill up a winning set of the game $\Ga'$. As for Breaker, due to the parity assumption, this strategy ensures that he gets all top vertices of $X_I$, which takes care of all winning sets that are not in the game $\Ga'$ by definition of $I$.
    \end{itemize}
    This ends the proof of Claim 1.
    \
    
    We define a \textit{good pattern} as a sequence of winning sets of the following form:
$$ (\{x_{i_1},x_{i_2}\}, \{x_{i_2},y_{i_3}\}, \{x_{i_3},y_{i_4}\}, \{x_{i_4},y_{i_5}\}, \ldots, \{x_{i_{\ell-1}},y_{i_{\ell}}\}), \,\,\text{with $\ell \geq 3$ and pairwise distinct $i_k$}. $$
	We define a \textit{winning pattern} as a good pattern (with notations as above) prolonged with a winning set of the form $\{x_{i_{\ell}},z\}$ where either: $z \in \{x_{i_2},\ldots,x_{i_{\ell-1}}\}$, $z \in \{y_{i_1},y_{i_2}\}$, or $z$ is a bottom vertex that does not appear in the good pattern. See Figure \ref{fig:WinningPattern}.
	

\begin{figure}[h]
    \centering
    \begin{tikzpicture}
        \node[v,minimum size =18 pt] (x1) at (0,0) {$x_{i_1}$};
        \node[v,minimum size =18 pt] (x2) at (1.25,0) {$x_{i_2}$};
        \node[v,minimum size =18 pt] (x3) at (2.5,0) {$x_{i_3}$};
        \node[v,minimum size =18 pt] (x4) at (3.75,0) {$x_{i_4}$};
        \node[v,minimum size =18 pt] (x5) at (5,0) {};
        \node[v,minimum size =18 pt] (x7) at (7.5,0) {};
        \node[v,minimum size =18 pt] (xl) at (8.75,0) {$x_{i_{\ell}}$};
        \node[v,minimum size =18 pt] (z) at (10,0) {$z$};
        
        \node[v,minimum size =18 pt] (y3) at (2.5,1.25) {$y_{i_3}$};
        \node[v,minimum size =18 pt] (y4) at (3.75,1.25) {$y_{i_4}$};
        \node[v,minimum size =18 pt] (y5) at (5,1.25) {};
        \node[v,minimum size =18 pt] (y7) at (7.5,1.25) {};
        \node[v,minimum size =18 pt] (yl) at (8.75,1.25) {$y_{i_{\ell}}$};

        \draw[->](x3)--(y3);
        \draw[->](x4)--(y4);
        \draw[->](x5)--(y5);
        \node at (6.25,0.75){...};
        \draw[->](x7)--(y7);
        \draw[->](xl)--(yl);

        \draw[thin, color = red] (0,-0.5) arc (270:90:0.5) -- (1.25,0.5) arc (90:-90:0.5) -- (0,-0.5);
        \draw[thin, color = red] (8.75,-0.5) arc (270:90:0.5) -- (10,0.5) arc (90:-90:0.5) -- (8.75,-0.5);
        \draw[thin, color = red] (1.53,-0.28) arc (315:135:0.4) -- (2.22,1.53) arc (135:-45:0.4)--cycle;
        \draw[thin, color = red] (2.78,-0.28) arc (315:135:0.4) -- (3.47,1.53) arc (135:-45:0.4)--cycle;
        \draw[thin, color = red] (4.03,-0.28) arc (315:135:0.4) -- (4.72,1.53) arc (135:-45:0.4)--cycle;
        \draw[thin, color = red] (7.78,-0.28) arc (315:135:0.4) -- (8.47,1.53) arc (135:-45:0.4)--cycle;
    \end{tikzpicture}
    \caption{A winning pattern where $z$ is a bottom vertex that does not appear in the good pattern.} 
    \label{fig:WinningPattern}
\end{figure}
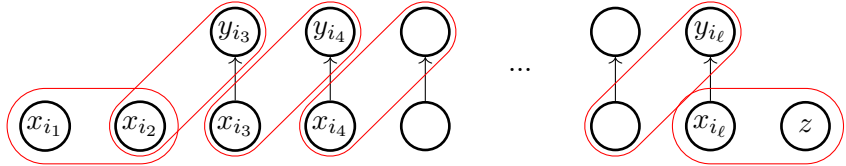
 
	$\quad$ Claim 2: If the game contains a winning pattern, then Maker wins.
	
	\ Let us prove Claim 2. Maker starts by claiming $x_{i_2}$, which forces Breaker to claim $x_{i_1}$ because of the winning set $\{x_{i_1},x_{i_2}\}$. Then Maker claims $x_{i_3}$, which forces Breaker to claim $y_{i_3}$ because of the winning set $\{x_{i_2},y_{i_3}\}$. Then Maker successively claims $x_{i_4},\ldots,x_{i_{\ell}}$, which forces Breaker to successively claim $y_{i_4},\ldots,y_{i_{\ell}}$ in similar fashion. If $z \in \{x_{i_2},\ldots,x_{i_{\ell-1}}\}$, then Maker has already filled up the winning set $\{x_{i_{\ell}},z\}$ and thus won the game. Otherwise, either $z \in \{y_{i_1},y_{i_2}\}$ or $z$ is a bottom vertex that does not appear in the good pattern: in both cases, $z$ is an available vertex, so Maker can now claim $z$ to win by filling up the winning set $\{x_{i_{\ell}},z\}$. This ends the proof of Claim 2.
	\
	
    $\quad$ Claim 3: If there exists at least one winning set and neither R1 nor R2 applies, then Maker wins.
    
    \
Let us prove Claim 3. Suppose there exists at least one winning set and neither R1 nor R2 applies. Let us first settle the case $t=0$ (all vertices are bottom): since R1 does not apply, there exists a winning set of the form $\{x_i\}$ or two winning sets of the form $\{x_i,x_j\}$ and $\{x_j,x_k\}$, so Maker wins in at most two rounds. Therefore, we now assume $r \geq 1$.
\\ We start by showing the existence (in all nontrivial cases) of a good pattern. There exists a winning set $W$ in which each vertex is a bottom vertex, otherwise R2 would apply to the set $I=\{1,\ldots,t\}$. Maker wins trivially if $|W|=1$, so assume $W=\{x_p,x_q\}$ for some $p \neq q$. Since R1 does not apply, we know that at least one of $x_p$ or $x_q$, say $x_q$, is part of a winning set $W' \neq W$. If $W'=\{x_q\}$, $W'=\{x_q,y_p\}$, $W'=\{x_q,y_q\}$ or $W'=\{x_q,x_r\}$ for some $r \neq p,q$, then Maker wins the game in at most two rounds, so assume $W'=\{x_q,y_r\}$ for some $r \neq p,q$. Therefore, there exist good patterns: $(W,W')$ is one.
\\ Now, we will show the existence (in all nontrivial cases) of a winning pattern. Let $I$ be the set of all indices of top vertices that appear in good patterns. Note that $I\neq \varnothing$ and $I \subseteq \{1,\ldots,t\}$. Since R2 does not apply to the set $I$, we know there exists $j \in I$ such that there exists a winning set of the form $\{x_j,z\}$ where $z \not\in \{y_i, i \in I\}$ (assuming $\{x_j\}$ is not a winning set, otherwise Maker trivially wins). Let $(\{x_{i_1},x_{i_2}\}, \{x_{i_2},y_{i_3}\}, \{x_{i_3},y_{i_4}\}, \{x_{i_4},y_{i_5}\}, \ldots, \{x_{i_{\ell-1}},y_{i_{\ell}}\})$ be a good pattern where $i_{\ell}=j$, which exists since $j \in I$. Since $z \not\in \{y_i, i \in I\}$, this good pattern does not contain $z$ as a top vertex, nor can it be prolonged into a longer good pattern by the winning set $\{x_j,z\}$. Moreover, we can also assume that $z \neq x_{i_1}$, otherwise Maker wins the game in two rounds by claiming $x_{i_1}$ as her first move. All in all, this leaves us with three possibilities: $z \in \{x_{i_2},\ldots,x_{i_{\ell-1}}\}$, $z \in \{y_{i_1},y_{i_2}\}$, or $z$ is a bottom vertex that does not appear in the good pattern. In all cases, we get a winning pattern by definition, so Claim 2 ensures that Maker wins. This ends the proof of Claim 3.
\

Finally, we should explain how to detect in $O(|X|^3)$ time  whether R2 applies or not. It suffices, for each index $1 \leq i_0 \leq t$, to answer the following question in $O(|X|^2)$ time: is there a set of indices containing $i_0$ and satisfying R2? Starting from the set $\{i_0\}$, iterate the operator $I \mapsto I \cup \bigcup_{i \in I} \{j \in \{1,\ldots,t\}\setminus I, \{x_i,y_j\} \in \WS \}$ until stability, and let $I_0$ be the resulting set. Note that $I_0$ can be computed in $O(|X|^2)$ time. Moreover, any set containing $i_0$ and satisfying R2 must include $I_0$. If $I_0$ satisfies R2, which can be checked in $O(|X|^2)$ time, then the answer is positive. Otherwise, it means there exists a winning set of the form $\{x_i\}$ or $\{x_i,x_j\}$ with $i \in I_0$, therefore the answer is negative.
\\ In conclusion, detecting whether R1 or R2 applies, and performing the reduction if so, can be done in $O(|X|^3)$ time. Obviously, at most $|X|$ such reductions can be performed in succession. When no reduction applies anymore, Claim 3 concludes: if at least one winning set remains, then Maker wins, otherwise Breaker trivially wins. \qedhere

 \end{itemize}
\end{proof}

\section{Conclusion and future work}

The current work introduces a new framework that opens the door to a variety of perspectives. Our analysis with respect to the parameters of the poset led to a first classification of the complexity of \popoproblem, which is a step towards a better understanding of the boundary between tractability and hardness. The results presented here directly induce a list of open problems that arise naturally in order to refine this boundary. Among them, we have identified the following two questions that seem the most relevant for us:

\begin{itemize}[noitemsep, topsep = 0cm]
\item When there is exactly one winning set of size $1$, the height of the poset makes a difference between tractability ($h=2$) and NP-hardness ($h=3$). Therefore, when $h=2$, it is worthwhile to examine how the conditions on the number $m$ of winning sets and their size $s$ impact the algorithmic complexity. A next study would be to investigate the case $m=2$ and $s=1$.
\item In the case of disjoint chains, an analysis of the complexity according to their width (i.e.~the number of chains) is a natural perspective. In particular, one may focus on the case $w=2$, even when restricted to $s\leq 3$ (since we have a polynomial-time algorithm for $s=2$). It would also be interesting to obtain a hardness results for disjoint chains.\end{itemize}

Going back to our initial motivation of giving a general framework for Connect-$k$ games, one could also examine the case of boards of even size for this game. For example, the famous case $(k,w,h)=(4,7,6)$, is known to be a Maker win in the Maker-Breaker convention since it a first player win in the Maker-Maker convention. However, a direct proof for the Maker-Breaker convention could be considered, that could also be extended to other sizes.

Finally, while we mostly studied the \popogames~in the Maker-Breaker convention, our definition can be transposed to all other conventions of positional games. It would therefore be interesting to perform a similar analysis for Maker-Maker or Avoider-Enforcer \popogames .

\bibliographystyle{plain}
\bibliography{references}

\begin{thebibliography}{10}

\bibitem{allen1990expert}
James~D. Allen.
\newblock Expert play in {C}onnect-{F}our.
\newblock 1990.
\newblock \url{https://tromp.github.io/c4.html}.

\bibitem{allen2010complete}
James~D. Allen.
\newblock {\em The Complete Book of Connect 4: History, Strategy, Puzzles}.
\newblock Puzzlewright Press, 2010.

\bibitem{allis1988knowledge}
Louis~Victor Allis.
\newblock {\em A Knowledge-Based Approach of {C}onnect-{F}our.}
\newblock Report IR-163. Vrije Universiteit Amsterdam, The Netherlands, 1988.

\bibitem{restrictedPosGames}
Pranav Avadhanam and Siddhartha~G. Jena.
\newblock Restricted positional games.
\newblock 2021.
\newblock Preprint, \url{https://doi.org/10.48550/arXiv.2108.12839}.

\bibitem{beck}
József Beck.
\newblock {\em Combinatorial Games: Tic-Tac-Toe Theory}.
\newblock Encyclopedia of Mathematics and its Applications. Cambridge
  University Press, 2008.

\bibitem{byskov2004maker}
Jesper~Makholm Byskov.
\newblock Maker-{M}aker and {M}aker-{B}reaker games are {PSPACE}-complete.
\newblock {\em BRICS Report Series}, 11(14), 2004.

\bibitem{Elkies1999OnNA}
Noam~D. Elkies.
\newblock {On numbers and endgames: combinatorial game theory in chess
  endgames}.
\newblock In {\em Games of No Chance, Proceedings of 7/94 MSRI Conference on
  Combinatorial Games}, pages 135--150. Cambridge University Press, 1996.

\bibitem{erdos}
Paul Erd\H{o}s and John~L. Selfridge.
\newblock {On a combinatorial game}.
\newblock {\em Journal of Combinatorial Theory}, 14:298--301, 1973.

\bibitem{Florian2}
Florian Galliot, Sylvain Gravier, and Isabelle Sivignon.
\newblock {Maker}-{Breaker} is solved in polynomial time on hypergraphs of rank
  3.
\newblock Preprint, \url{https://arxiv.org/abs/2209.12819}, 2022.

\bibitem{Hales1963}
Alfred~W. Hales and Robert~I. Jewett.
\newblock Regularity and positional games.
\newblock {\em Trans. Am. Math. Soc}, 106:222--229, 1963.

\bibitem{positionalgames}
Dan Hefetz, Michael Krivelevich, Milo\v{s} Stojakovi\'{c}, and Tibor Szab\'{o}.
\newblock {\em {Positional Games}}, volume~44 of {\em Oberwolfach Seminars}.
\newblock Birkhäuser (Springer), 2014.

\bibitem{Karp1972}
Richard~M. Karp.
\newblock Reducibility among combinatorial problems.
\newblock In R.E. Miller, J.W. Thatcher, and J.D. Bohlinger, editors, {\em
  Complexity of Computer Computations: Proceedings of a symposium on the
  Complexity of Computer Computations, March 20--22, 1972}, chapter~9, pages
  85--103. Springer US, Boston, MA, 2004.

\bibitem{rahman2021}
Md~Lutfar Rahman and Thomas Watson.
\newblock {6-uniform {M}aker-Breaker game is {PSPACE}-complete}.
\newblock In {\em Proceedings of the 38th International Symposium on
  Theoretical Aspects of Computer Science ({STACS 2021})}, volume 187 of {\em
  Leibniz International Proceedings in Informatics (LIPIcs)}, pages 57:1--15,
  2021.

\bibitem{schaefer}
Thomas~J. Schaefer.
\newblock {On the complexity of some two-person perfect-information games}.
\newblock {\em Journal of Computer and System Sciences}, 16:185--225, 1978.

\bibitem{stockmeyer1973}
Larry~J Stockmeyer and Albert~R Meyer.
\newblock Word problems requiring exponential time (preliminary report).
\newblock In {\em Proceedings of the Fifth Annual ACM Symposium on Theory of
  Computing}, pages 1--9, 1973.

\bibitem{tromp2015}
John Tromp.
\newblock {John's Connect Four Playground}.
\newblock \url{http://tromp.github.io/c4/c4.html}.

\bibitem{tromp2008solving}
John Tromp.
\newblock Solving {C}onnect-4 on medium board sizes.
\newblock {\em ICGA Journal}, 31(2):110--112, 2008.

\end{thebibliography}

\newpage

\section*{Appendix}

\subsection{Outcome of the disjoint union of two games}

\begin{theorem}\label{thm:union}
    The possible outcomes for the disjoint union of two \popogames~is given in Figure~\ref{table: union}.
\end{theorem}

\begin{figure}[!htb]
\begin{center}
\begin{tabular}{ c || c | c | c | c }
 $even \backslash even$ & $\OM$ & $\ON$ & $\OP$ & $\OB$ \\ 
 \hline \hline
 $\OM$ & $\OM$ & $\OM$        & $\OM$   & $\OM$ \\ \hline  
 $\ON$ & $\OM$ & $\OM, \ON$   & $\OM$   & $\OM, \ON$ \\ \hline 
 $\OP$ & $\OM$ & $\OM$        & $\OP$   & $\OP$ \\ \hline
 $\OB$ & $\OM$ & $\OM, \ON$   & $\OP$   & $\OP, \OB$  \\
\end{tabular}
\end{center}

\vspace{.3cm}
\begin{center}
\begin{tabular}{ c || c | c | c | c }
 $even \backslash odd$ & $\OM$ & $\ON$ & $\OP$ & $\OB$ \\ 
 \hline \hline
 $\OM$ & $\OM$ & $\OM, \ON$         & $\OM$         & $\OM, \ON$ \\  \hline
 $\ON$ & $\OM$ & $\OM, \ON$         & $\OM$         & $\OA$ \\ \hline
 $\OP$ & $\OM$ & $\ON$              & $\OM$         & $\ON$ \\ \hline
 $\OB$ & $\OM$ & $\ON$              & $\OM, \OP$    & $\ON, \OB$  \\
\end{tabular}
\end{center}

\vspace{.3cm}
\begin{center}
\begin{tabular}{ c || c | c | c | c }
 $odd \backslash odd$ & $\OM$ & $\ON$ & $\OP$ & $\OB$ \\ 
 \hline \hline
 $\OM$ & $\OM$      & $\OM$              & $\OM$         & $\OM, \ON$ \\  \hline
 $\ON$ & $\OM$      & $\OM, \OP$         & $\OM, \ON$    & $\OA$ \\ \hline
 $\OP$ & $\OM$      & $\OM, \ON$         & $\OM$         & $\OM, \ON$ \\ \hline
 $\OB$ & $\OM, \ON$ & $\OA$              & $\OM, \ON$    & $\OA$  \\
\end{tabular}
\end{center}
    \caption{Possible outcomes for a disjoint union of two \popogames~depending on the parity of the number of vertices in both games. We write "$\OA$" when all four outcomes are possible.}
    \label{table: union}
\end{figure}

For convenience, a \popogame~will be deemed \emph{even} (resp. \emph{odd}) if it has an even (resp. odd) number of vertices. We start by showing that the outcomes which do not appear in Figure~\ref{table: union} are indeed impossible.

\begin{lemma}\label{lem:union_cases}
    Let $\Ga$ and $\Ga'$ be two \popogames. The following assertions hold:
    \begin{enumerate}
        \item If $o(\Ga) \in \{\OM,\OP\}$ and $o(\Ga') \in \{\OM,\OP\}$, then $o(\Ga \cup \Ga') \in \{\OM,\OP\}$.
        \item If $o(\Ga) \in \{\OM,\OP\}$ and $o(\Ga') \in \{\OM,\ON\}$, then $o(\Ga \cup \Ga') \in \{\OM,\ON\}$.
        \item If $\Ga$ is even and $o(\Ga') \in \{\OM,\ON\}$, then $o(\Ga \cup \Ga') \in \{\OM,\ON\}$.
        \item If $\Ga$ is even and $o(\Ga') \in \{\OM,\OP\}$, then $o(\Ga \cup \Ga') \in \{\OM,\OP\}$.
        \item If $\Ga$ is odd and $o(\Ga') \in \{\OM,\OP\}$, then $o(\Ga \cup \Ga') \in \{\OM,\ON\}$.
        \item If both $\Ga$ and $\Ga'$ are odd and have an outcome in $\{\OM,\ON\}$, then $o(\Ga \cup \Ga') \in \{\OM,\OP\}$.
        \item If both $\Ga$ and $\Ga'$ are even and have an outcome in $\{\OB,\OP\}$, then $o(\Ga \cup \Ga') \in \{\OB,\OP\}$.
        \item If $\Ga$ is even, $\Ga'$ is odd, $o(\Ga) \in \{\OB,\OP\}$ and $o(\Ga') \in \{\OB,\ON\}$, then $o(\Ga \cup \Ga') \in \{\OB,\ON\}$.
    \end{enumerate}
\end{lemma}

\begin{proof}
    First of all, let us establish four important claims.
    \
    
    $\quad$ Claim 1: Let $\Ga_1$ and $ \Ga_2$ be two \popogames. If Maker has a winning strategy in both $\Ga_1$ and $\Ga_2$ as second player, then Maker has a winning strategy in $\Ga_1 \cup \Ga_2$ as second player.
    \
    
    $\quad$ Claim 2: Let $\Ga_1$ and $ \Ga_2$ be two \popogames. If $\Ga_1$ and $\Ga_2$ are even and Breaker has a winning strategy in both $\Ga_1$ and $\Ga_2$ as second player, then Breaker has a winning strategy in $\Ga_1 \cup \Ga_2$ as second player.
    \
        
    $\quad$ Claim 3: Let $\Ga_1$ and $ \Ga_2$ be two \popogames. If $\Ga_1$ is even and Maker has a winning strategy in $\Ga_2$ as first player, then Maker has a winning strategy in $\Ga_1 \cup \Ga_2$ as first player.
     \
    
    $\quad$ Claim 4: Let $\Ga_1$ and $ \Ga_2$ be two \popogames. If $\Ga_1$ is even and Maker has a winning strategy in $\Ga_2$ as second player, then Maker has a winning strategy in $\Ga_1 \cup \Ga_2$ as second player.
    \
    
    Let us prove those claims. Claim 1 is rather straightforward. Maker always plays in the same component Breaker has just played in, according to her respective winning strategy. By the time one of the two components has all its vertices claimed, Maker will necessarily have filled a winning set in that component. Claim 2 is similar: Breaker always plays in the same component Maker has just played in, according to his respective winning strategy. The parity ensures that the moves alternate correctly in both components (even after one of them has all its vertices claimed), so that Breaker successfully acts as second player in both components, thus avoiding filling any winning set. Finally, Claims 3 and 4 are direct consequences of Lemma \ref{lem:emptyUnion}: indeed, Maker would have a winning strategy even if $\Ga_1$ had no winning set, so the same strategy also works if $\Ga_1$ has winning sets. We can now go on with the proof of the lemma:

    \begin{enumerate}
        \item Since Maker has a winning strategy in both $\Ga$ and $\Ga'$ as second player, she has a winning strategy in $\Ga \cup \Ga'$ as second player by Claim 1, which means $o(\Ga \cup \Ga') \in \{\OM,\OP\}$.
        \item Since Maker has a winning strategy in $\Ga'$ as first player, there exists a vertex $u$ in $\Ga'$ such that, in the game $\Ga''$ that results from $\Ga'$ after Maker claims $u$, Maker has a winning strategy as second player. If Maker claims $u$ as her first move in $\Ga \cup \Ga'$, we obtain the game $\Ga \cup \Ga''$, in which Maker has a winning strategy as second player by Claim 1. Therefore, Maker has a winning strategy in $\Ga \cup \Ga'$ as first player, which means $o(\Ga \cup \Ga') \in \{\OM,\ON\}$.
        \item Since $\Ga$ is even and Maker has a winning strategy in $\Ga'$ as first player, she has a winning strategy in $\Ga \cup \Ga'$ as first player by Claim 3, which means $o(\Ga \cup \Ga') \in \{\OM,\ON\}$.
        \item Since $\Ga$ is even and Maker has a winning strategy in $\Ga'$ as second player, she has a winning strategy in $\Ga \cup \Ga'$ as second player by Claim 4, which means $o(\Ga \cup \Ga') \in \{\OM,\OP\}$.
        \item Playing first, Maker may claim an arbitrary vertex in $\Ga$: we get a disjoint union $\Ga'' \cup \Ga'$ where $\Ga''$ is even and Maker has a winning strategy in $\Ga'$ as second player. By Claim 4, Maker has a winning strategy in $\Ga'' \cup \Ga'$ as second player, so she has a winning strategy in $\Ga \cup \Ga'$ as first player, which means $o(\Ga \cup \Ga') \in \{\OM,\ON\}$.
        \item Assume by symmetry that, playing first, Breaker claims some vertex in $\Ga$: we get a disjoint union $\Ga'' \cup \Ga'$ where $\Ga''$ is even and Maker has a winning strategy in $\Ga'$ as first player. By Claim 3, Maker has a winning strategy in $\Ga'' \cup \Ga'$ as first player, so she has a winning strategy in $\Ga \cup \Ga'$ as second player, which means $o(\Ga \cup \Ga') \in \{\OM,\OP\}$.
        \item Since $\Ga$ and $\Ga'$ are even and Breaker has a winning strategy in both $\Ga$ and $\Ga'$ as second player, he has a winning strategy in $\Ga \cup \Ga'$ as second player by Claim 2, which means $o(\Ga \cup \Ga') \in \{\OB,\OP\}$.
        \item Since Breaker has a winning strategy in $\Ga'$ as first player, there exists a vertex $u$ in $\Ga'$ such that, in the even game $\Ga''$ that results from $\Ga'$ after Breaker claims $u$, Breaker has a winning strategy as second player. If Breaker claims $u$ as his first move in $\Ga \cup \Ga'$, we obtain the game $\Ga \cup \Ga''$, in which Breaker has a winning strategy as second player by Claim 2. Therefore, Breaker has a winning strategy in $\Ga \cup \Ga'$ as first player, which means $o(\Ga \cup \Ga') \in \{\OB,\ON\}$. \qedhere
    \end{enumerate}
\end{proof}

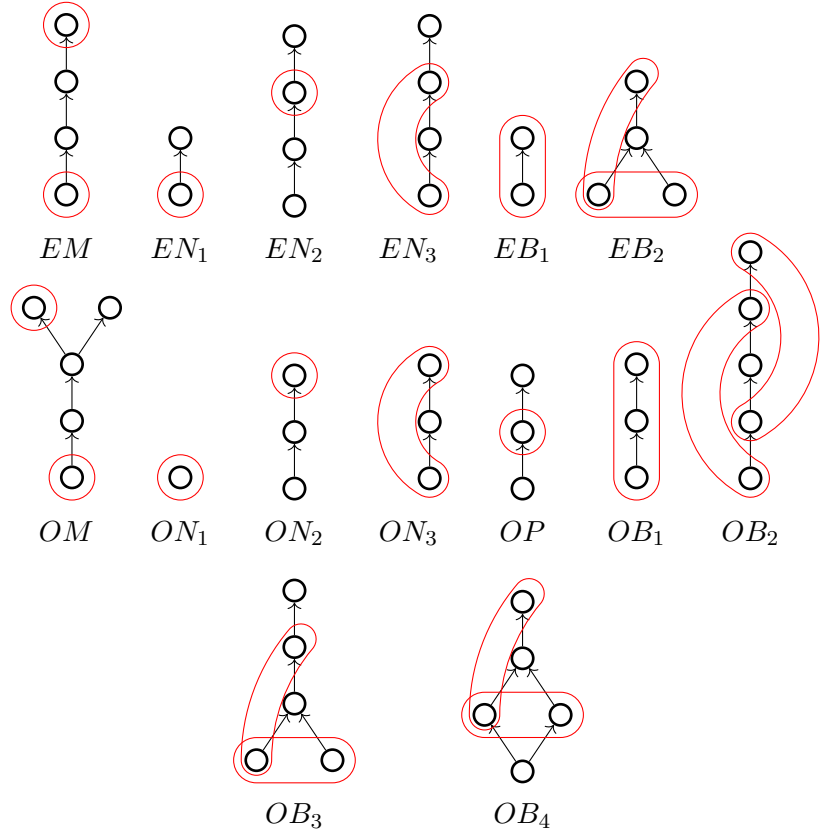
\begin{figure}[h]
    \centering
    \begin{tikzpicture}
        \node[anchor = south, label= below:$EM$] at (-1.5,0){
            \begin{tikzpicture}
                \node[circle, inner sep=0, minimum size =8 pt, line width = 1pt, draw=black, fill=white, anchor = center] (u1) at (0,0){};
                \node[circle, inner sep=0, minimum size =8 pt, line width = 1pt, draw=black, fill=white,, anchor = center] (u2) at (0,0.75){};
                \node[circle, inner sep=0, minimum size =8 pt, line width = 1pt, draw=black, fill=white,, anchor = center] (u3) at (0,1.5){};
                \node[circle, inner sep=0, minimum size =8 pt, line width = 1pt, draw=black, fill=white,, anchor = center] (u4) at (0,2.25){};
                \draw[->] (u1) -- (u2);
                \draw[->] (u2) -- (u3);
                \draw[->] (u3) -- (u4);
                \draw[thin, color = red] (0,0) circle (0.3);
                \draw[thin, color = red] (0,2.25) circle (0.3);
            \end{tikzpicture}
        };
        \node[anchor = south, label= below:$EN_1$] at (0,0){
            \begin{tikzpicture}
                \node[circle, inner sep=0, minimum size =8 pt, line width = 1pt, draw=black, fill=white,, anchor = center] (u1) at (0,0){};
                \node[circle, inner sep=0, minimum size =8 pt, line width = 1pt, draw=black, fill=white,, anchor = center] (u2) at (0,0.75){};
                \draw[->] (u1) -- (u2);
                \draw[thin, color = red] (0,0) circle (0.3);
            \end{tikzpicture}
        };
        \node[anchor = south, label= below:$EN_2$] at (1.5,0){
            \begin{tikzpicture}
                \node[circle, inner sep=0, minimum size =8 pt, line width = 1pt, draw=black, fill=white,, anchor = center] (u1) at (0,0){};
                \node[circle, inner sep=0, minimum size =8 pt, line width = 1pt, draw=black, fill=white,, anchor = center] (u2) at (0,0.75){};
                \node[circle, inner sep=0, minimum size =8 pt, line width = 1pt, draw=black, fill=white,, anchor = center] (u3) at (0,1.5){};
                \node[circle, inner sep=0, minimum size =8 pt, line width = 1pt, draw=black, fill=white,, anchor = center] (u4) at (0,2.25){};
                \draw[->] (u1) -- (u2);
                \draw[->] (u2) -- (u3);
                \draw[->] (u3) -- (u4);
                \draw[thin, color = red] (0,1.5) circle (0.3);
            \end{tikzpicture}
        };
        \node[anchor = south, label= below:$EN_3$] at (3,0){
            \begin{tikzpicture}
                \node[circle, inner sep=0, minimum size =8 pt, line width = 1pt, draw=black, fill=white,, anchor = center] (u1) at (0,0){};
                \node[circle, inner sep=0, minimum size =8 pt, line width = 1pt, draw=black, fill=white,, anchor = center] (u2) at (0,0.75){};
                \node[circle, inner sep=0, minimum size =8 pt, line width = 1pt, draw=black, fill=white,, anchor = center] (u3) at (0,1.5){};
                \node[circle, inner sep=0, minimum size =8 pt, line width = 1pt, draw=black, fill=white,, anchor = center] (u4) at (0,2.25){};
                \draw[->] (u1) -- (u2);
                \draw[->] (u2) -- (u3);
                \draw[->] (u3) -- (u4);
                \draw [thin, color = red] (0.125,0.216) arc (60:-120:0.25) arc (240:120:1.116) arc (120:-60:0.25) arc (120:240:0.616);
            \end{tikzpicture}
        };
        \node[anchor = south, label= below:$EB_1$] at (4.5,0){
            \begin{tikzpicture}
                \node[circle, inner sep=0, minimum size =8 pt, line width = 1pt, draw=black, fill=white,, anchor = center] (u1) at (0,0){};
                \node[circle, inner sep=0, minimum size =8 pt, line width = 1pt, draw=black, fill=white,, anchor = center] (u2) at (0,0.75){};
                \draw[->] (u1) -- (u2);
                \draw[thin, color = red] (0.3,0) arc (0:-180:0.3) -- (-0.3,0.75) arc (180:0:0.3)--cycle;
            \end{tikzpicture}
        };
        \node[anchor = south, label= below:$EB_2$] at (6,0){
            \begin{tikzpicture}
                \node[circle, inner sep=0, minimum size =8 pt, line width = 1pt, draw=black, fill=white,, anchor = center] (u1) at (0,0){};
                \node[circle, inner sep=0, minimum size =8 pt, line width = 1pt, draw=black, fill=white,, anchor = center] (u2) at (1,0){};
                \node[circle, inner sep=0, minimum size =8 pt, line width = 1pt, draw=black, fill=white,, anchor = center] (u3) at (0.5,0.75){};
                \node[circle, inner sep=0, minimum size =8 pt, line width = 1pt, draw=black, fill=white,, anchor = center] (u4) at (0.5,1.5){};
                \draw[->] (u1) -- (u3);
                \draw[->] (u2) -- (u3);
                \draw[->] (u3) -- (u4);
                \draw[thin, color = red] (0,0.3) arc (90:270:0.3) -- (1,-0.3) arc (270:450:0.3)--cycle;
                \draw[thin, color = red] (0.2,0) arc (0:-180:0.2) arc (180:140:2.7) arc (140:-40:0.2) arc (140:180:2.3);
            \end{tikzpicture}
        };
        \node[anchor = south, label= below:$OM$] at (-1.5,-3.75){
            \begin{tikzpicture}
                \node[circle, inner sep=0, minimum size =8 pt, line width = 1pt, draw=black, fill=white,, anchor = center] (u1) at (0,0){};
                \node[circle, inner sep=0, minimum size =8 pt, line width = 1pt, draw=black, fill=white,, anchor = center] (u2) at (0,0.75){};
                \node[circle, inner sep=0, minimum size =8 pt, line width = 1pt, draw=black, fill=white,, anchor = center] (u3) at (0,1.5){};
                \node[circle, inner sep=0, minimum size =8 pt, line width = 1pt, draw=black, fill=white,, anchor = center] (u4) at (0.5,2.25){};
                \node[circle, inner sep=0, minimum size =8 pt, line width = 1pt, draw=black, fill=white,, anchor = center] (u5) at (-0.5,2.25){};
                \draw[->] (u1) -- (u2);
                \draw[->] (u2) -- (u3);
                \draw[->] (u3) -- (u4);
                \draw[->] (u3) -- (u5);
                \draw[thin, color = red] (0,0) circle (0.3);
                \draw[thin, color = red] (u5) circle (0.3);
            \end{tikzpicture}
        };
        \node[anchor = south, label= below:$ON_1$] at (0,-3.75){
            \begin{tikzpicture}
                \node[circle, inner sep=0, minimum size =8 pt, line width = 1pt, draw=black, fill=white,, anchor = center] (u1) at (0,0){};
                \draw[thin, color = red] (0,0) circle (0.3);
            \end{tikzpicture}
        };
        \node[anchor = south, label= below:$ON_2$] at (1.5,-3.75){
            \begin{tikzpicture}
                \node[circle, inner sep=0, minimum size =8 pt, line width = 1pt, draw=black, fill=white,, anchor = center] (u1) at (0,0){};
                \node[circle, inner sep=0, minimum size =8 pt, line width = 1pt, draw=black, fill=white,, anchor = center] (u2) at (0,0.75){};
                \node[circle, inner sep=0, minimum size =8 pt, line width = 1pt, draw=black, fill=white,, anchor = center] (u3) at (0,1.5){};
                \draw[->] (u1) -- (u2);
                \draw[->] (u2) -- (u3);
                \draw[thin, color = red] (0,1.5) circle (0.3);
            \end{tikzpicture}
        };
        \node[anchor = south, label= below:$ON_3$] at (3,-3.75){
            \begin{tikzpicture}
                \node[circle, inner sep=0, minimum size =8 pt, line width = 1pt, draw=black, fill=white,, anchor = center] (u1) at (0,0){};
                \node[circle, inner sep=0, minimum size =8 pt, line width = 1pt, draw=black, fill=white,, anchor = center] (u2) at (0,0.75){};
                \node[circle, inner sep=0, minimum size =8 pt, line width = 1pt, draw=black, fill=white,, anchor = center] (u3) at (0,1.5){};
                \draw[->] (u1) -- (u2);
                \draw[->] (u2) -- (u3);
                \draw [thin, color = red] (0.125,0.216) arc (60:-120:0.25) arc (240:120:1.116) arc (120:-60:0.25) arc (120:240:0.616);
            \end{tikzpicture}
        };
        \node[anchor = south, label= below:$OP$] at (4.5,-3.75){
            \begin{tikzpicture}
                \node[circle, inner sep=0, minimum size =8 pt, line width = 1pt, draw=black, fill=white,, anchor = center] (u1) at (0,0){};
                \node[circle, inner sep=0, minimum size =8 pt, line width = 1pt, draw=black, fill=white,, anchor = center] (u2) at (0,0.75){};
                \node[circle, inner sep=0, minimum size =8 pt, line width = 1pt, draw=black, fill=white,, anchor = center] (u3) at (0,1.5){};
                \draw[->] (u1) -- (u2);
                \draw[->] (u2) -- (u3);
                \draw[thin, color = red] (u2) circle (0.3);
            \end{tikzpicture}
        };
        \node[anchor = south, label= below:$OB_1$] at (6,-3.75){
            \begin{tikzpicture}
                \node[circle, inner sep=0, minimum size =8 pt, line width = 1pt, draw=black, fill=white,, anchor = center] (u1) at (0,0){};
                \node[circle, inner sep=0, minimum size =8 pt, line width = 1pt, draw=black, fill=white,, anchor = center] (u2) at (0,0.75){};
                \node[circle, inner sep=0, minimum size =8 pt, line width = 1pt, draw=black, fill=white,, anchor = center] (u3) at (0,1.5){};
                \draw[->] (u1) -- (u2);
                \draw[->] (u2) -- (u3);
                \draw[thin, color = red] (0.3,0) arc (0:-180:0.3) -- (-0.3,1.5) arc (180:0:0.3)--cycle;
            \end{tikzpicture}
        };
        \node[anchor = south, label= below:$OB_2$] at (7.5,-3.75){
            \begin{tikzpicture}
                \node[circle, inner sep=0, minimum size =8 pt, line width = 1pt, draw=black, fill=white,, anchor = center] (u1) at (0,0){};
                \node[circle, inner sep=0, minimum size =8 pt, line width = 1pt, draw=black, fill=white,, anchor = center] (u2) at (0,0.75){};
                \node[circle, inner sep=0, minimum size =8 pt, line width = 1pt, draw=black, fill=white,, anchor = center] (u3) at (0,1.5){};
                \node[circle, inner sep=0, minimum size =8 pt, line width = 1pt, draw=black, fill=white,, anchor = center] (u4) at (0,2.25){};
                \node[circle, inner sep=0, minimum size =8 pt, line width = 1pt, draw=black, fill=white,, anchor = center] (u5) at (0,3){};
                \draw[->] (u1) -- (u2);
                \draw[->] (u2) -- (u3);
                \draw[->] (u3) -- (u4);
                \draw[->] (u4) -- (u5);
                \draw [thin, color = red] (0.125,0.216) arc (60:-120:0.25) arc (240:120:1.55) arc (120:-60:0.25) arc (120:240:1.05);
                \draw [thin, color = red] (-0.125,0.966) arc (120:300:0.25) arc (-60:60:1.55) arc (60:240:0.25) arc (60:-60:1.05);
            \end{tikzpicture}
        };
        \node[anchor = south, label= below:$OB_3$] at (1.5,-7.5){
            \begin{tikzpicture}
                \node[circle, inner sep=0, minimum size =8 pt, line width = 1pt, draw=black, fill=white,, anchor = center] (u1) at (0,0){};
                \node[circle, inner sep=0, minimum size =8 pt, line width = 1pt, draw=black, fill=white,, anchor = center] (u2) at (1,0){};
                \node[circle, inner sep=0, minimum size =8 pt, line width = 1pt, draw=black, fill=white,, anchor = center] (u3) at (0.5,0.75){};
                \node[circle, inner sep=0, minimum size =8 pt, line width = 1pt, draw=black, fill=white,, anchor = center] (u4) at (0.5,1.5){};
                \node[circle, inner sep=0, minimum size =8 pt, line width = 1pt, draw=black, fill=white,, anchor = center] (u5) at (0.5,2.25){};
                \draw[->] (u1) -- (u3);
                \draw[->] (u2) -- (u3);
                \draw[->] (u3) -- (u4);
                \draw[->] (u4) -- (u5);
                \draw[thin, color = red] (0,0.3) arc (90:270:0.3) -- (1,-0.3) arc (270:450:0.3)--cycle;
                \draw[thin, color = red] (0.2,0) arc (0:-180:0.2) arc (180:140:2.7) arc (140:-40:0.2) arc (140:180:2.3);
            \end{tikzpicture}
        };
        \node[anchor = south, label= below:$OB_4$] at (4.5,-7.5){
            \begin{tikzpicture}
                \node[circle, inner sep=0, minimum size =8 pt, line width = 1pt, draw=black, fill=white,, anchor = center] (u0) at (0.5,-0.75){};
                \node[circle, inner sep=0, minimum size =8 pt, line width = 1pt, draw=black, fill=white,, anchor = center] (u1) at (0,0){};
                \node[circle, inner sep=0, minimum size =8 pt, line width = 1pt, draw=black, fill=white,, anchor = center] (u2) at (1,0){};
                \node[circle, inner sep=0, minimum size =8 pt, line width = 1pt, draw=black, fill=white,, anchor = center] (u3) at (0.5,0.75){};
                \node[circle, inner sep=0, minimum size =8 pt, line width = 1pt, draw=black, fill=white,, anchor = center] (u4) at (0.5,1.5){};
                \draw[->] (u0) -- (u1);
                \draw[->] (u0) -- (u2);
                \draw[->] (u1) -- (u3);
                \draw[->] (u2) -- (u3);
                \draw[->] (u3) -- (u4);
                \draw[thin, color = red] (0,0.3) arc (90:270:0.3) -- (1,-0.3) arc (270:450:0.3)--cycle;
                \draw[thin, color = red] (0.2,0) arc (0:-180:0.2) arc (180:140:2.7) arc (140:-40:0.2) arc (140:180:2.3);
            \end{tikzpicture}
        };
    \end{tikzpicture}
    \caption{A selection of \popogames~that illustrate all possibilities for the outcome of a disjoint union.}
    \label{fig:union}
\end{figure}

\begin{proof}[Proof of Theorem~\ref{thm:union}]
    For the disjoint union, by Lemma~\ref{lem:union_cases}, any outcome that does not appear in Figure~\ref{table: union} is impossible. Conversely, we now show that all the outcomes listed in Figure~\ref{table: union} can be achieved. Consider the \popogames~described in Figure~\ref{fig:union}: $EM$, $EN_1$, $EN_2$, $EN_3$, $EB_1$ and $EB_2$ are even games, while $OM$, $ON_1$, $ON_2$, $ON_3$, $OP$, $OB_1$, $OB_2$, $OB_3$ and $OB_4$ are odd games. Moreover, the reader can verify that we have $o(EM) = o(OM)= \OM$, $o(EN_1) =o(EN_2) =o(EN_3) = o(ON_1)= o(ON_2)= o(ON_3)= \ON$, $o(OP)= \OP$ and $o(EB_1) =o(EB_2) =o(OB_1) = o(OB_2)= o(OB_3)= o(OB_4)= \OB$. Now, consider the following disjoint unions:

    \begin{itemize}
        \item We have $o(EN_1 \cup EN_1) = \OM$ and $o(EN_1 \cup EN_2) = \ON$, therefore it is possible for the disjoint union of two even games of outcome $\ON$ to have outcome $\OM$ or $\ON$.
        
        \item We have $o(EN_1 \cup EB_1) = \ON$ and $o(EN_1 \cup EB_2) = \OM$, therefore it is possible for the disjoint union of two even games of outcome $\ON$ and $\OB$ respectively to have outcome $\OM$ or $\ON$.
        
        \item We have $o(EB_1 \cup EB_1) = \OB$ and $o(EB_2 \cup EB_2) = \OP$, therefore it is possible for the disjoint union of two even games of outcome $\OB$ to have outcome $\OB$ or $\OP$.  

        \item We have $o(EM \cup ON_1) = \OM$ and $o(EM \cup ON_2) = \ON$, therefore it is possible for the disjoint union of an even game of outcome $\OM$ and an odd game of outcome $\ON$ to have outcome $\OM$ or $\ON$.
        
        \item We have $o(EN_1 \cup ON_1) = \OM$ and $o(EN_1 \cup ON_2) = \ON$, therefore it is possible for the disjoint union of an even game of outcome $\ON$ and an odd game of outcome $\ON$ to have outcome $\OM$ or $\ON$. 
        
        \item We have $o(EB_2 \cup OP) = \OM$ and $o(EB_1 \cup OP) = \OP$, therefore it is possible for the disjoint union of an even game of outcome $\OB$ and an odd game of outcome $\OP$ to have outcome $\OM$ or $\OP$.  
        
        \item We have $o(EM \cup OB_3) = \OM$ and $o(EM \cup OB_1) = \ON$, therefore it is possible for the disjoint union of an even game of outcome $\OM$ and an odd game of outcome $\OB$ to have outcome $\OM$ or $\ON$.  
        
        \item We have $o(EN_1 \cup OB_3) = \OM$, $o(EN_1 \cup OB_1) = \ON$, $o(EN_2 \cup OB_1) = \OP$  and $o(EN_3 \cup OB_1) = \OB$, therefore it is possible for the disjoint union of an even game of outcome $\ON$ and an odd game of outcome $\OB$ to have any of the four outcomes.
        
        \item We have $o(EB_2 \cup OB_1) = \ON$ and $o(EB_1 \cup OB_1) = \OB$, therefore it is possible for the disjoint union of an even game of outcome $\OB$ and an odd game of outcome $\OB$ to have outcome $\ON$ or $\OB$. 
        
        \item We have $o(OM \cup OB_3) = \OM$ and $o(OM \cup OB_1) = \ON$, therefore it is possible for the disjoint union of two odd games of outcome $\OM$ and $\OB$ respectively to have outcome $\OM$ or $\ON$.
        
        \item We have $o(ON_1 \cup ON_1) = \OM$ and $o(ON_2 \cup ON_2) = \OP$, therefore it is possible for the disjoint union of two odd games of outcome $\ON$ to have outcome $\OM$ or $\OP$.
        
        \item We have $o(ON_2 \cup OP) = \OM$ and $o(ON_1 \cup OP) = \ON$, therefore it is possible for the disjoint union of two odd games of outcome $\ON$ and $\OP$ respectively to have outcome $\OM$ or $\ON$.
        
        \item We have $o(ON_2 \cup OB_2) = \OM$, $o(ON_1 \cup OB_1) = \ON$, $o(ON_2 \cup OB_1) = \OP$  and $o(ON_3 \cup OB_1) = \OB$, therefore it is possible for the disjoint union of two odd games of outcome $\ON$ and $\OB$ respectively to have any of the four outcomes.
        
        \item We have $o(OP \cup OB_2) = \OM$ and $o(OP \cup OB_1) = \ON$, therefore it is possible for the disjoint union of two odd games of outcome $\OP$ and $\OB$ respectively to have outcome $\OM$ or $\ON$.
        
        \item We have $o(OB_2 \cup OB_2) = \OM$, $o(OB_1 \cup OB_3) = \ON$, $o(OB_4 \cup OB_4) = \OP$  and $o(OB_1 \cup OB_1) = \OB$, therefore it is possible for the disjoint union of two odd games of outcome $\OB$ to have any of the four outcomes. \qedhere
        
    \end{itemize}
\end{proof}

\subsection{Posets of height 2, winning sets of size 1 in Maker-Maker convention}

We establish {\sf PSPACE}-completeness, in the Maker-Maker convention, for a family of \popogames~which Theorem \ref{thm:topverticesWS} has shown to be solved in polynomial time in the Maker-Breaker convention.

\begin{theorem}\label{thm:MMsize1depth2}
For \popogames~in the Maker-Maker convention, deciding whether the first (resp. second) player has a winning strategy is {\sf PSPACE}-complete, even when restricted to instances where: the poset is of height 2, all the winning sets are of size 1, and all the nonminimal elements are winning sets.
\end{theorem}

\begin{proof}
We perform a reduction from {\sc Avoid True}, which is defined as follows. Given a $2$-POSDNF formula~$\phi$ on a set of variables $V$, two players, Alice and Bob, alternately choose a variable in $V$ (which has not been chosen yet) and set it to True, until one of the clauses of $\phi$ is satisfied. The player who first turns $\phi$ to True loses. Note that no draw is possible. Deciding the outcome of {\sc Avoid True} has been proved to be {\sf PSPACE}-complete by Schaefer \cite{schaefer}.

Let $(V, \phi)$ be an instance of {\sc Avoid True}. Denote by $x_1, \dots, x_n$ the variables in $V$ and by $C_1, \dots, C_m$ the clauses of $\phi$. We build a \popogame~$\Ga = (X,\WS,\poset) $ as follows:

\begin{itemize}
    \item For each variable $x_i \in V$, we add a vertex $v_i$ in $X$.
    \item For each clause $C_j \in \phi$, we add a vertex $u_j$ in $X$ and the winning set $\{u_j\}$ in $\WS$.
    \item For each variable $x_i \in C_j$, we add $v_i < u_j$ in $\poset$.
\end{itemize}

We prove that the first player in {\sc Avoid True} has a winning strategy in $(V, \phi)$ if and only if the first player in $\Ga$ has a winning strategy (which is equivalent to the second player in $\Ga$ not having a winning strategy, since no draw will be possible in $\Ga$). 

We will name the players Alice and Bob in both games. By symmetry, suppose Alice has a winning strategy $\strat$ in {\sc Avoid True}, going first (resp. second). We define a strategy for her in $\Ga$ as follows. Alice considers that any vertex $v_i$ being claimed corresponds to the variable $x_i$ being set to True in {\sc Avoid True}. Whenever it is Alice's turn, she claims the vertex $v_i$ corresponding to the variable $x_i$ that she would have set to True according to $\strat$. Since $\strat$ is a winning strategy for Alice, Bob will be the first in {\sc Avoid True} to set some variable $x_i$ to True which turns a clause $C_j$ to True. This means that $v_i$ is the last predecessor of $u_j$ in $\Ga$. Therefore, Alice simply claims $u_j$ to win.
\end{proof}

\end{document}